\documentclass{ws-ijitdm}
\usepackage{tabularx} 
\usepackage{float}
\usepackage{hyperref}
\usepackage{graphicx} 
\makeatletter
\def\footnoterule{\kern-3\p@
  \hrule \@width 2in \kern 2.6\p@} 
\makeatother
\begin{document}

\title{Algorithms for Multi-Criteria Decision-Making and Efficiency Analysis Problems}

\author{Fuh-Hwa Franklin Liu \footnotemark{}\footnotetext{Corresponding author; fliu@nycu.edu.tw} }
\address{Professor Emeritus, Department of Industrial Engineering and Management,  National Yang-Ming Chiao Tung University, Hsin Chu City, Taiwan 300, Republic of China}
\author{Su-Chuan Shih \footnotemark{}\footnotetext{scshih@gm.pu.edu.tw}}
\address{Associate Professor, Department of Business Administration, Providence University, \\ Tai Chung City, Taiwan 433, Republic of China}

\maketitle
\begin{history}
\received{7 June 2024}
\end{history}
\begin{abstract} 
Multi-criteria decision-making (MCDM) problems involve the evaluation of alternatives based on various minimization and maximization criteria. Similarly, efficiency evaluation (EA) methods assess decision-making units (DMUs) by analyzing their input consumption and output production. MCDM and EA methods face challenges in managing alternatives and DMUs with varying capacities across different criteria (inputs and outputs). That leads to performance assessments often skewed by subjective biases in criteria weighting. We introduce two innovative scenarios utilizing linear programming-based Virtual Gap Analysis (VGA) models to address these limitations. This dual-scenario approach aims to mitigate traditional biases, offering robust solutions for comprehensively assessing alternatives and DMUs. Our methodology allows for the influential ranking of alternatives in MCDM problems and enables each DMU to adjust its input and output ratios to achieve efficiency.
\end{abstract}
\section{Research Background and Objectives} \label{sec:1}
\subsection{The Background} \label{sec:1.1}
Decision-makers across various sectors are frequently confronted with complex multi-criteria decision-making (MCDM) scenarios, necessitating the evaluation of alternatives against a broad range of criteria. Miettinen surveyed methods to visualize alternatives in MCDM problems\cite{1}. Sahoo and Goswami\cite{2} comprehensively reviewed MCDM methodologies, highlighting significant challenges posed by the subjective assessment of alternatives due to the heterogeneous capacities of criteria. The efficiency analysis (EA) method, Data Envelopment Analysis (DEA) is a primary tool employed to address MCDM challenges\cite{2}. Opricović and Tzeng\cite{3} compared the DEA model and the MCDM VIKOR method. The terms "alternative" and "decision-making unit" (DMU) are used interchangeably.  Conflicting criteria in EA and MCDM methods involve input and output volumes that must be minimized and maximized. Linear programming-based DEA models evaluate each DMU's inefficiency score against its peers, as detailed by Cooper et al.\cite{4}, the DMU under evaluation is referred to as DMU-o.

\indent Based on a robust theoretical foundation for over fifty years, DEA models\cite{5} often oversimplify evaluations by assuming homogeneity among DMUs. Each model employs a set of artificial goal weights to estimate DMU-o's criteria weights regardless of input and output configurations, potentially leading to inaccuracies in operational contexts that reflect the diverse realities of DMUs. Appendix A provides a concise overview of basic DEA models, distinguishing them from VGA models. These distinctions are further summarized in Table \ref{table:4} in Appendix B.

\indent Despite these limitations, Data Envelopment Analysis (DEA) models have inspired the development of innovative Virtual Gap Analysis (VGA) scenarios. VGA, grounded in linear programming theories and operating without assumptions, provides comprehensive solutions for assessing each DMU-o. These solutions include:
\begin{romanlist}
\item Adjustment Ratio and Virtual Unit Price: VGA calculates the adjustment ratio and virtual unit price for each input and output, offering detailed insights into resource allocation.
\item Intensities of DMUs: The intensities of DMUs are analyzed to understand their operational dynamics.
\item Reference Peers: VGA identifies reference peers for each DMU-o, facilitating comparative analysis and benchmarking.
\item Efficiency Score: An efficiency score is assigned to each DMU-o, reflecting its performance relative to the best practices observed within the dataset.
\item Sensitivity analysis: VGA allows for sensitivity analysis on the data utilized in the linear models, enabling an examination of how changes in data affect the outcomes. This feature is precious for ensuring the robustness and reliability of the results.
\end{romanlist}

\indent We anticipate that DMU-o will find these solutions achievable in its realistic practice, enhancing the practical applicability of the analysis. In contrast, traditional statistics-based EA methods and existing MCDM methods do not offer such a comprehensive and detailed analysis, highlighting the advantages of the VGA approach.
\subsection{Our Algorithms to Solve MCDM and EA Problems} \label{sec:1.2}
We have developed two scenarios to address MCDM and EA problems. In the first scenario, we propose the first and second VGA models to estimate the inefficiency score of each DMU-o, excluding and including the sum of intensities condition (SIC), which equates the total variable intensities of DMUs to an SIC scalar. We developed a four-phase procedure, where the first three phases measure the range of the SIC scalar. In the final phase, DMU-o selects the SIC scalar to determine achievable adjustment ratios for inputs and outputs.

\indent Each inefficient DMU has a realistic efficiency score between 0 and 1, alongside reduction and expansion ratios for inputs and outputs. Each efficient DMU has zero adjustment ratios and an efficiency score of 1. In the second scenario, we assess each efficient DMU to estimate its super-efficiency, which is significantly greater than 1. The super-efficiency converges to 1 when the efficient DMU adjusts inputs and outputs' estimated expansion and reduction ratios.

\indent The third and fourth VGA models employed in the second scenario are modifications of the first two VGA models. Whether inefficient or efficient, each DMU comprehensively measures the achievable adjustment ratio and realistic virtual price per unit of each input and output.

\indent The MCDM decision-making group (DMG) utilizes these two scenarios to select alternatives with a distinct approach. In the second scenario, the DMG selects a final SIC scalar for each efficient alternative according to the four-phase procedure. Under specific subjective considerations in decision-making, the MCDM DMG employs this systematic method to compare alternatives effectively.
\subsection{Objectives of the Research}\label{sec:1.3}
This research introduces innovative linear programming-based VGA models under two distinct scenarios for evaluating DMU-o. The VGA models are designed to accommodate the heterogeneity of DMUs, thereby promising a comprehensive evaluation of DMU-o performance aligned with the systematic \textit{virtual unit price} of each criterion. The specific objectives of this study are:
\begin{romanlist}[(ii)]
\item Ensure that DMU-o aligns with reference peers of similar input and output configurations.
\item Enable DMU-o to mirror reference peers to precisely calibrate the estimations of inputs and outputs' virtual unit prices and \textit{adjustment ratios}.
\item Allow DMU-o to select the final SIC \textit{scalar} within the model to calibrate all adjustment ratios (AARs) achievable in practice.
\item For the first scenario, achieve the following outcomes:
\begin{itemize}
\item Assign an inefficiency score to DMU-o realistically ranging between 0 and 1.
\item Express criteria weaknesses through estimated adjustment ratios, identifying the weakest criterion as the input or output with the largest adjustment ratio.
\end{itemize}
\item Advance from the first to the second scenario to evaluate each identified efficient DMU, DMU-o, leading to:
\begin{itemize}
\item Assign a super-efficiency score to DMU-o that is realistically larger than 1.
\item Observe that as DMU-o adjusts the estimated expansion and reduction ratios of inputs and outputs, its estimated super-efficiency score deteriorates to 1.
\item Express criteria strengths through estimated adjustment ratios, identifying the strongest criterion as the input or output with the largest adjustment ratio.
\end{itemize}
\item Enable decision-makers in MCDM to effectively compare efficient alternatives based on their criteria strengths and super-efficiencies.
\end{romanlist}
 
 \indent Section \ref{sec:1.2} introduces the substantial strengths of the first VGA scenario. The first scenario employs the first and second VGA models, accomplishing objectives (i)–(iv). These two models are further illustrated in Sections \ref{sec:2.2} and \ref{sec:2.3}. The third and fourth VGA models are integrated into the second scenario to fulfill objectives (i)–(iii) and (v)–(vi). Details of these two VGA models are offered in Section \ref{sec:3.1} and \ref{sec:3.2}. 
\subsection{Early VGA Developments}\label{sec:1.4}
Hwang\cite{6} contributed to establishing CRS and VRS VGA models, where the former excludes and includes the SIC, equating to 1. However, issues arose concerning the virtual goal price for DMU-o and the relationship between these models. The virtual goal price was subjectively assigned. Past research, including works by Liu and Liu\cite{7,8}, applied the CRS VGA model in two-phase production systems. The authors have consulted on unpublished theses that utilized these models. 
\subsection{Paper Organization} \label{sec:1.5}
This paper is organized as follows: Section 2 introduces the first and second VGA models for assessing inefficiency. Section 3 discusses the third and fourth VGA models for determining super-efficiency. Section 4 describes the four-phase VGA procedure. Section 5 presents a numerical example to illustrate the two VGA scenarios. Finally, Section 6 concludes the paper and discusses the findings. Appendix A points out the reasons why DEA models have incomplete evaluations. Appendix B compares DEA and VGA models. Appendix C introduces the computational tools of VGA scenarios.
\section {First Scenario: Evaluating DMUs' Inefficiencies} \label{sec:2} The first scenario assesses the inefficiencies of DMUs by the first and second VGA models in a four-phase procedure.
\subsection {Notations} \label{sec:2.1}
 The set  $J$,  $I$, and  $R$, respectively, denotes the  $n$ DMUs,  $m$ inputs, and  $s$ outputs, where $n \ge2(m+s)$. The \textit {decision matrix}, denoted as  $(\textit{X}$, $\textit{Y})$, consists of column vectors $(x_j,y_j)$, in which the elements  $(x_{ij},y_{rj})$ represent positive \textit {observed volume} of each \textit {input-i} and \textit {output-r} of each DMU-j.

  \indent The \textit{adjustment ratio} of each input-i and output-r are decision variables  $Q_{io}$ and  $P_{ro}$. In evaluating the inefficient DMU-o by the VGA models in Section \ref{sec:2}, input-i reduces $Q_{io}$ and output-r expands $P_{ro}$.   $\pi_{jo}$ denotes the variable intensity of DMU-j. The variable virtual unit price for input-i and output-r are  $(v_{io}, u_{ro})$. Symbols such as   $Q_o$, $P_o$, $v_o$, $u_o$, and $\pi_o$  represent the vectors' decision variables. 
  
  \indent The symbol  $w_o^c$ is the \textit{scalar unit price} corresponding to the SIC.  $\kappa_o^c$ is the SIC scalar of the TSc model.  $\omega_o^a$ is the \textit{virtual scalar price}, $\kappa_o^c \times w_o^c$. 
  
  \indent We developed a \textit{two-step process} to determine a systematic \textit{unified goal price},  $\tau_o$, measured in the virtual currency \$. The superscripts  \# and   $\star$ indicate the optimal solutions for a decision variable in the VGA model for Step I and II.
  
  \indent The symbol of a decision variable with superscript  "$\star$" represents its optimal solution.

 \subsection{Pure Technical Efficiency (PT) Model-the First VGA Model} \label{sec:2.2} This model measures the adjustment ratio and virtual unit price of each input and output of DMU-o alongside the reference peers' intensities.

\textbf{The dual program} of the PT model aims to measure the maximum total adjustment price (TAP) of each DMU-o.	
\begin{equation} \label{Eq:1}
	\delta_o^{PT\star} = \max_{Q_o, P_o, \pi} \sum_{\forall i\in I} Q_{io} \tau_o  + \sum_{\forall r\in R}P_{ro} \tau_o, \forall o \in J;
 \end{equation}
\vspace{-1 em}
 \begin{equation} \label{Eq:2}
	s.t. \sum_{\forall j\in J} x_{ij} \pi_{jo} + Q_{io} x_{io} =x_{io}  , \forall i\in I;
 \end{equation}
\vspace{- 1 em}
\begin{equation} 
 - \sum_{\forall j \in J} y_{rj} \pi_{jo} + P_{ro} y_{ro} = - y_{ro}, \forall r \in R;	
\label{Eq:3}
 \end{equation}
 \vspace{- 1.5 em}
\begin{equation} \label{Eq:4}
\pi_o, Q_o  , P_o  \ge 0.	
 \end{equation}
Each \textit{adjustment condition} for each input (or output), Eqs. (\ref{Eq:2}) and (\ref{Eq:3}), ensures that an input reduces (or output expands)  $Q_{io}$ (or $P_{ro}$) to the benchmark,  $\sum_{\forall j\in J} x_{ij} \pi_{jo}$ (or $\sum_{\forall j\in J} y_{rj} \pi_{jo})$. Each adjustment condition mirrors the primal program's variable virtual price per input-i (or output-r) unit, $v_{io}$ (or $u_{ro}$). 
Eq. (\ref{Eq:1}) seeks to maximize the sum of variable AARs multiplied by the specified goal price. Distinguished by positive intensities, reference peers determine the AARs under these adjustment conditions. The estimated AARs are realistic since DMU-o and its reference peers maintain similar input and output configurations.
 
\indent \textbf{The primal program} of the PT model aims to measure the total virtual gap (TVG):
\begin{equation}
	\Delta_o^{PT\star} = \min_{v_o, u_o} \sum_{\forall i\in I} v_{io} x_{io}  - \sum_{\forall r\in R} u_{ro} y_{ro} , \forall o \in J;	
 \label{Eq:5}
 \end{equation}
 \vspace{- 1.5 em}
 \begin{equation}
	s.t. \sum_{\forall i\in I} v_{io} x_{ij}  - \sum_{\forall r\in R}u_{ro}  y_{rj}  \ge 0, \forall j\in J; \label{Eq:6}
 \end{equation}
 \vspace{- 1.5 em}
 \begin{equation}
	x_{io} v_{io} \ge\tau_o, \forall i\in I;	
 \label{Eq:7}
 \end{equation}
 \vspace{- 1.5 em}
\begin{equation}\label{Eq:8}
	y_{ro} u_{ro} \ge \tau_o, \forall r \in R;
 \end{equation}
 \vspace{- 1.5 em}
\begin{equation}
	v_o, u_o \quad free.	
 \label{Eq:9}
 \end{equation}

 \indent Each \textit{virtual price condition} in Eq. (\ref{Eq:7}) and (\ref{Eq:8}) restricts the \textit{virtual price} of input-i and output-r, $x_{io} v_{io}$ and $y_{ro} u_{ro}$, to a unified lower bound \textit{goal price} $\tau_o$,  denominated in \textit{virtual currency} (\$). 
 Each \textit{virtual price condition}  correlates with the dual program's variable adjustment ratio, $Q_{io}$ (or $P_{ro}$).
 
 \indent The \textit{virtual gap condition} (\ref{Eq:6}) limits each DMU-j to a total virtual gap, the excess of \textit{virtual input} and \textit{virtual output}, $\sum_{\forall i\in I} v_{io} x_{ij}$ - $\sum_{\forall r\in R}u_{ro}  y_{rj}$, to a minimum of \$0. This condition aligns with the dual program's dimensionless variable intensity, $\pi_{jo}$. Efficient reference peers exhibit a zero virtual gap.  

\indent  \textbf{DMU-o determines its unified goal price in two steps}: (1) The initial goal price is set at \$1, and (2) the second goal price is derived from the optimal solutions. Using the second goal price, DMU-o is expected to achieve an estimated virtual input of \$1 and ensure that its reference peers share similar criteria configurations.
  
\indent In Step I and Step II , substitute $\tau_o$  by  $\tau_o^\#$ = \(\$1\) and  $\tau_o^\star$   = $\$ \bar{t}$, respectively. An inefficient DMU-o will obtain the optimal solutions $\$ 0 \leq  \delta_o^{PT\#}=\Delta_o^{PT\#} $, $\$ 0 \leq \delta_o^{PT\star }=\Delta_o^{PT\star } \leq \$1$,  and $(Q_o^\star ,P_o^\star , \pi^\star ) = (Q_o^\#, P_o^\#,\pi^\# )$. Since $\delta_o^{PT\#} : \delta_o^{PT\star }=\tau_o^\#:\tau_o^\star   = \$ 1 : \$ \bar{t}$, therefore, $\delta_o^{PT\star }$  is a function of $\bar{t}$. Similarly, in  (\ref{Eq:5}), $v_o  x_o$  is a function of $\bar{t}$. Determine the dimensionless value of $\bar{t}$ to have $v_o^\star x_o $=$\$$1 so that $\delta_o^{PT\star }$  is ensured between \$0 and \$1. According to $\tau_o^\# :\tau_o^\star   = \$ 1 : \$ \bar{t}$, and yields $v_{o}^\# x_o : v_o^\star  x_o =\$ 1 : \$ \bar{t}$. Use (\ref{Eq:10}) to obtain the dimensionless value of  $\bar{t}$. The solutions of Step II  are $v_o^\star   x_o =\$1 \textrm{  and  }u_o^\star   y_o \leq \$1$.
\begin{equation}
	\bar{t} = \$ 1/v_{o}^\#  x_o \textrm {  and  } \tau_o^\star   = \$ \bar{t}. 
	\label{Eq:10}
  \end{equation}
The optimal solution of (\ref{Eq:1}) as depicted in (\ref{Eq:11}), the total adjustment price (TAP) of DMU-o in Step II equals $\delta_{xo}^{PT\star}$ plus $\delta_{yo}^{PT\star}$. (\ref{Eq:12}) expresses the observed multiple inputs and outputs with distinct measurement units are aggregated into the $[\textit{pure virtual input  (pvInput)},  \textit{pure virtual output (pvOutput)}]$ of DMU-j in Steps I and II, $(\alpha_j^\#, \beta_j^\#)$ and $(\alpha_j^\star, \beta_j^\star)$.  The minimized  \textit{total virtual gap (TVG)} of DMU-o equals to the pvInput minus pvOutput,  $(\alpha_o^\# -\beta_o^\#)$ and  $(\alpha_o^\star -\beta_o^\star  )$ in Steps I and-II .   
 \begin{equation}
 \begin{aligned}
&\$0 \leq \delta_o^{PT\#} = \sum_{\forall i\in I} Q_{io} ^\# \times \$1  + \sum_{\forall r\in R}P_{ro} ^\#  \times \$1 ; \\ 
&\$0 \leq \delta_o^{PT\star } = \sum_{\forall i\in I} Q_{io}^\star  \tau_o ^\star  + \sum_{\forall r\in R} P_{ro}^\star  \tau_o ^\star =\delta_{xo}^{PT\star} + \delta_{yo}^{PT\star}\leq \$1.
  \end{aligned}  \label{Eq:11} 
 \end{equation} 
 
\indent \textbf{Aggregate inputs and outputs into pvInput and pvOutput.}  
 \begin{equation} \label {Eq:12}
\begin{aligned}
\$0 \leq (TVG)=(pvInput)-(pvOutput) ;\\
\$0 \leq  \Delta_j^{PT\#}  = v_o^\# x_j - u_o^\# y_j = \alpha_j^\# - \beta_j^\# ; \\
\$0 \leq  \Delta_j^{PT\star }  = v_o^\star  x_j - u_o^\star  y_j =\alpha_j^\star  - \beta_j^\star \leq \$1.
\end{aligned}
 \end{equation}

\indent The primal program's objective Eq. (\ref{Eq:5}) is to minimize the virtual gap of DMU-o, ranging between \$0 and \$1. At the same time, the dimensionless inefficiency score is calculated as the ratio of the estimated virtual gap to the estimated virtual input, ideally ranging between 0 and 1.
\begin{theorem}\label{thm:1}
{The PT model fulfills Condition C1, the estimated PT efficiency between 0 and 1.}
\end{theorem}
\noindent\textbf{Proof.} The PT model estimates the technical parameters $v_o^\star$ and $u_o^\star$ to have the minimum virtual gap  $\Delta_o^{PT\star }$. Which is converted into the $\textit{pure technical}$ (PT) inefficiency    $F_o^{PT\star }$ according to (\ref{Eq:13}), where  $E_o^{PT\star }$is the efficiency. 
\begin{equation}
\begin{aligned}\label{Eq:13} 
	0\le E_o^{PT\star } = 1 -F_o^{PT\star } = \beta_o^\star  /\alpha_o^\star  \le 1; F_o^{PT\star } = \Delta_o^{PT\star }/\alpha_o^\star = (\alpha_o^\star -\beta_o^\star)/\alpha_o^\star.	  
\end{aligned}
\end{equation}
$E_o^{PT\star }$ and $F_o^{PT\star }$ are between zero and one,  where $\alpha_o^\star = \$1$ and $\beta_o^\star  \leq \$1$. \hfill $\square$

\indent Certainly, normalizing Step I solutions to attain Step II  solutions is a common practice in mathematical optimization problems.  (\ref{Eq:14}) likely represents this normalization process, where the Step I solutions are adjusted or transformed to achieve the solutions in Step II. This normalization might involve scaling or modifying Step I solutions to fulfill specific conditions required for Step II.
\begin{equation} \begin{aligned}
  (\delta_o^{PT\star },  \Delta_o^{PT\star },  v_o^\star, u_o^\star) = \bar{t} \times (\delta_o^{PT\#}, \Delta_o^{PT\#}, v_o^\#, u_o^\#).
\label {Eq:14} \end{aligned}
\end{equation}
The symbol  $\mathcal{E}_o^{PT}$  denotes the set of \textit{reference DMUs} of DMUo in the PT evaluations. Each reference DMU-j is efficient, and its estimated intensity, $\pi_{jo}^\star >0$. At the same time, the estimated intensity of the other inefficient DMUs, $\pi_{jo}^\star =0$. DMU-o estimates the benchmark of each input-i and output-r,  $\widehat{x}_{io}^{PT\star}$ and  $\widehat{y}_{ro}^{PT\star}$ via  (\ref{Eq:15}). DMU-o mimics reference DMUs with their estimated intensities,  $\pi_{jo}^\star$, $\forall j\in \mathcal{E}_o^{PT}$.
\begin{equation}
\begin{aligned}
	&\widehat{x}_{io}^{PT\star} = \sum_{\forall j\in \mathcal{E}_o^{PT}} x_{ij}  \pi_{jo}^\star  = x_{io} (1 - Q_{io}^\star ), \forall i\in I; \\
	&\widehat{y}_{ro} ^{PT\star} = \sum_{\forall j\in \mathcal{E}_o^{PT}}y_{rj} \pi_{jo}^\star =y_{ro} (1 + P_{ro}^\star  )
 \forall r\in R.
\end{aligned}
\label {Eq:15}
\end{equation}
\begin{theorem}\label{thm:2}
{ PT model fulfills Condition C2, the adjusted PT efficiency equals 1.}
\end{theorem}
\noindent \textbf{Proof.} Assessing DMU-o with the adjusted input-i and output-r, \(\widehat{x}_{io}^{PT\star}\) and \(\widehat{y}_{ro} ^{PT\star} \), as shown in  (\ref{Eq:15}), will have the PT efficiency   $\widehat{E}_o^{PT\star }$ equals 1.  \hfill $\square$
Let the total of estimated intensities be the first SIC scalar, $\kappa_o^1$, computed via  (\ref{Eq:16}). 
\begin{equation}
    \sum_{\forall j\in \mathcal{E}_o^{PT}}\pi_{jo}^\star   = \kappa_o^1 . 
 \label {Eq:16} \end{equation} 
\subsection {Technical and Scalar Choice (TSc) Model} \label{sec:2.3}
According to  (\ref{Eq:15}), the benchmarks of inputs and outputs are the function of the intensities of DMUs. One can manipulate the SIC scalar to choose the proper assessments. Adding (\ref{Eq:20}) to the PT model's dual program will have the TSc model's dual program. The SIC, $\sum_{\forall j\in J} \pi_{jo} =\kappa_o^c$, corresponds to the free-in-sign decision variable, $w_o^c$. The optimal solutions are affected by the SIC. The virtual goal price $\tau_o^\star $ is systematically determined in two steps.

\noindent TAP (dual) program of the TSc model:	
\begin{equation}\label {Eq:17} 
	\delta_o^{TSc\star}= \max_{Q_o, P_o} \sum_{\forall i\in I} Q_{io} \tau_o + \sum_{\forall r\in R}P_{ro}   \tau_o, \forall o \in J;	
\end{equation}
 \vspace{- 1 em}
 \begin{equation} \label{Eq:18}
	s.t. \sum_{\forall j\in J } x_{ij} \pi_{jo} + Q_{io} x_{io}  = x_{io}, \forall i \in I;	
 \end{equation}
  \vspace{-1 em}
\begin{equation}\label {Eq:19}
	- \sum_{\forall j\in J} y_{rj} \pi_{jo} + P_{ro} y_{ro} = - y_{ro}  , \forall r \in R;	
\end{equation}
 \vspace{-1 em}
 \begin{equation}\label {Eq:20}
	 \sum_{\forall j\in J } \pi_{jo}  = \kappa_o^c;
  \end{equation}
   \vspace{-1.5em}
  \begin{equation}\label {Eq:21}
	\pi_o,  Q_o, P_o  \ge 0.	
\end{equation}
  TVG (primal) program of the TSc model:
\begin{equation}\begin{aligned}
	\Delta_o^{TSc\star} =  \min_{x_o, v_o, w_o^c }  \sum_{\forall i\in I} v_{io} x_{io}  &- \sum_{\forall r\in R} u_{ro} y_{ro} + \kappa_o^c w_o^c, \forall o \in J;
 \label {Eq:22}
 \end{aligned}
 \end{equation}
  \vspace{-1.5 em}
 \begin{equation}\label {Eq:23}
s.t. \sum_{\forall i\in I} v_{io} x_{ij} - \sum_{\forall r\in R} u_{ro} y_{rj}  + 1 w_o^c  \ge 0, \forall j \in J; 
\end{equation} 
 \vspace{-1.5 em}
 \begin{equation}\label{Eq:24}
x_{io} v_{io} \ge\tau_o,\forall i \in I; 
\end{equation}
 \vspace{-1.5em}
 \begin{equation}\label {Eq:25}
y_{ro} u_{ro} \ge \tau_o, \forall r \in R;	
 \end{equation}
  \vspace{-1.5em}
 \begin{equation}\label {Eq:26}
v_o, u_o, w_o^c \quad free.
\end{equation}
\indent In Steps I and II , substitute \(\tau_o\) by \(\tau_o^\# = \$ 1\) and \(\tau_o^\star   = \$ \bar{t}\) will solve the TAP program comprehensively with the solutions \((Q_o^\star , P_o^\star , \pi^\star )=(Q_o^\#, P_o^\#, \pi^\# )\). In  (\ref{Eq:22}), the elements $(\kappa_o^c  w_o^\#)$ and \((\kappa_o^c  w_o^\star) \) are repressed by the symbols  \(\omega_o^{c\#}\) 
and  \(\omega_o^{c\star}\), the \textit{ virtual scalar \$  (vScalar)} in Steps I and II. 

\indent As depicted in  (\ref{Eq:27}),  \(\gamma_o\) and \((1- \gamma_o)\) are the proportions of the total adjustment prices of inputs and outputs. If \(w_o^{c\#}\) (or \(w_o^{c\star }\)) equals \(\$ 0\), let $\gamma_o$ equal 0.5. Partition the \textit{vScalar} into two parts to reflect the effects of the SIC on inputs and outputs of DMU-o. We denote them as \emph{vScalar in inputs (ivScalar)} and \emph{vScalar in outputs (ovScalar)}; vScalar equals ivScalar plus ovScalar.
\begin{equation}\begin{aligned}
	\gamma_o : (1-\gamma_o) = \tau_o^\# \times\sum_{\forall i\in I}Q_{io}^\# :\tau_o^\# \times\sum_{\forall r\in R}P_{ro}^\# = \tau_o^\star  \times\sum_{\forall i\in I}Q_{io}^\star  :\tau_o^\star \times\sum_{\forall r\in R}P_{ro}^\star .
 \end{aligned}\label {Eq:27}
 \end{equation}
\indent The minimized TVG of DMU-o in Steps II of (\ref{Eq:22}) is expressed as (\ref{Eq:28}); Step I has a similar expression. The estimated TVG, $\Delta_o^{TSc\star}$ equals the \emph{gray virtual input (gvInput} )  minus \emph{gray virtual output (gvOutput} ) plus vScalar.  $\Delta_o^{TSc\star}$ equals the  \textit{affected virtual Input (avInput)} minus \textit{affected virtual Output   (avOutput)}.  
 \begin{equation}\begin{aligned}\label {Eq:28}
\$0 \leq \Delta_o^{TSc\star} &= v_o^\star x_o-u_o^\star y_o+\omega_o^{c\star } = [v_o^\star x_o+(1-\gamma_o)\omega_o^{c\star }]-(u_o^\star   y_o-\gamma_o\omega_o^{c\star }) \\
&=[gvInput+ivScalar]-[gvOutput-ovScalar]\\
&= avInput-avOutput  = \alpha_o^{c\star} - \beta_o^{c\star}  \leq \$1; \\
\$0 \leq \Delta_o^{TSc\#} &=\alpha_o^{c\#} - \beta_o^{c\#}.
 \end{aligned} \end{equation}

 \indent Similarly, the estimated \textit{vScalar} of Step II  in  (\ref{Eq:23}), \(1w_o^*\), is decomposed into two components; Step I solutions has the similar expressions. As shown in  (\ref{Eq:29}), the TVG of each DMU-j is expressed as the gap between the \textit{avInput} and \textit{avOutput}.
  \begin{equation}\begin{aligned} \label{Eq:29}
\$0 \leq \Delta_j^{TSc\star} &= v_o^\star x_j-u_o^\star y_j+ 1w_o^{c\star }\\
&= [v_o^\star x_j+(1-\gamma_o) 1w_o^{c\star }]-(u_o^\star   y_j-\gamma_o 1w_o^{c\star })
= \alpha_j^{c\star } - \beta_j^{c\star }, \forall j \in J.
 \end{aligned} \end{equation}
 Each DMU-j is expressed by the pair of \emph{(avInput, avOutput)}, the virtual scales. The symbol  \(\mathcal{E}_o^{TS}\) denotes the set of reference peers in evaluating DMU-o by the TSc model, where each reference DMU-j has \(\pi_{jo}^\star >0\). The other inefficient DMU-j has \(\pi_{jo}^\star =0\). In  (\ref{Eq:29}), any best DMU-j belongs to \(\mathcal{E}_o^{TS}\) has  $\alpha_j^{c\star } = \beta_j^{c\star }$, while the other remaining DMUs have  $\alpha_j^{c\star } > \beta_j^{c\star }$. 
\begin{theorem}\label{thm:3} {The TSc model fulfills Condition C1, the estimated TSc efficiency score between 0 and 1.}
\end{theorem}
\noindent\textbf{Proof.} Using (\ref{Eq:30}) would have the dimensionless value of $\bar{t}$ by setting the \textit{avInput} in Step II $\alpha_o^{c\star}$ equals $\$1$. Therefore, the TVG in Step II , $\Delta_o^{TSc\star }$, should be between \$0 and \$1.
\begin{equation} \begin{aligned} \label {Eq:30}
	\$\bar{t }= \$ 1 / \alpha_o^{c\#}= \$ 1 / [v_o^\#  x_o +(1 - \gamma_o)\omega_o^{c\#}] =\tau_o^\star.
\end{aligned} \end{equation}
Because $\tau_o^\#:\tau_o^\star  = \$1 : \$\bar{t}$, therefore, if $\tau_o^\# $ = $\$\bar{t}$, then \(\alpha_o^{c\star}\) equals \$1. According to (\ref{Eq:28}),  $\alpha_o^{c\star}$ -  $\beta_o^{c\star}$ is greater than \$0. Using (\ref{Eq:31}), the estimated virtual gap  $\Delta_o^{TSc\star }$ is converted into the maximized TSc inefficiency score,  $E_o^{TSc\star}$, should be between 0 and 1, in which the sum of   \(F_o^{TSc\star }\) and \(E_o^{TSc\star }\) equals 1. 
\begin{equation}
\begin{aligned}
	0 < F_o^{TSc\star } &=1-E_o^{TSc\star } = \Delta_o^{TSc\star} / \alpha_o^{c\star}=(\alpha_o^{c\star}  - \beta_o^{c\star}) /\alpha_o^{c\star}  < 1;\\
 0 < E_o^{TSc\star } &=\frac{
	 u_o^\star y_o-\gamma_o \omega_o^{c\star}}{v_{o}^\star x_o+(1-\gamma_o)\omega_o^{c\star} } = \beta_o^{c\star}/\alpha_o^{c\star} <1.	\end{aligned}
 \label {Eq:31}
  \end{equation}
 \hfill $\square$ 

\indent Normalizing the solutions of Step I by the dimensionless value $\bar{t}$ would obtain the solutions of Step II, as shown in  (\ref{Eq:32}). 
\begin{equation}\begin{aligned}\label {Eq:32}
	(\Delta_o^{TSc\star }, \delta_o^{TSc\star }, v_o^\star , u_o^\star , w_o^{c\star })
	= \quad\bar{t}\times(\Delta_o^{TSc\#},\delta_o^{TSc\#} ,v_o^\#, u_o^\# , w_o^{c\#}).
\end{aligned}
\end{equation}
\indent Use  (\ref{Eq:33}) to compute the benchmark of each performance index,  \(\widehat{x}_{io}^{TSc\star}\) , and  \(\widehat{y}_{ro}^{TSc\star}\). DMU-o imitates the best peers with their estimated intensities \((\pi_{jo}^\star )\). DMU-o would become efficient with the benchmarks. 
\begin{equation}
\begin{aligned}\label {Eq:33} 
	\widehat{x}_{io}^{TSc\star} &= \sum_{\forall j\in \mathcal{E}_o^{TS}}x_{ij}\pi_{jo}^\star = x_{io}(1 - Q_{io}^\star ), \forall i\in I;\\
	\widehat{y}_{ro} ^{TSc\star} &= \sum_{\forall j\in \mathcal{E}_o^{TS}}y_{rj}\pi_{jo}^\star = y_{ro}(1 + P_{ro} ^\star ), \forall r\in R .   \end{aligned} \end{equation}
\begin{theorem}\label{thm:4}{The TSc model fulfills Condition C2, the adjusted DMU-o has the TSc efficiency equals 1.}
\end{theorem}
\noindent\textbf{Proof.} Assessing DMU-o with the adjusted input-i and output-r as depicted in  (\ref{Eq:33}), \(\widehat{x}_{io}^{TSc\star}\) and \(\widehat{y}_{ro} ^{TSc\star} \), by the TSc model will have $\delta_o^{TSc\star}=\$0$ and the TSc efficiency   $\widehat{E}_o^{TSc\star }=1$. \hfill $\square$
\subsection{Duality Properties}\label{sec:2.4}
We formulated VGA models based on the theories of linear programming (Dantzig and Thapa, 1997). When the optimal TVG and TAP of the two VGA models are equal, as stated in (\ref{Eq:34}), it suggests an equilibrium or balance between the defined conditions or objectives within these models. This equality could signify an alignment between the virtual gap and slack price considerations, indicating a harmonized solution or relationship between the gap and price metrics.
\begin{equation}\begin{aligned} \label{Eq:34}
\delta_o^{PT\star}=\Delta_o^{PT\star} \quad \textrm {and}   \quad
\delta_o^{TSc\star}=\Delta_o^{TSc\star}.
  \end{aligned} \end{equation}
\indent The strong complementary slackness conditions of the two VGA models are shown in (\ref{Eq:35}) and (\ref{Eq:36}).
\begin{equation} \begin{aligned} \label {Eq:35}
\relax [\sum_{\forall j\in \mathcal{E}_o^{VGA}} x_{ij} \pi_{jo}^\star - x_{io} (1 - Q_{io}^\star)] v_{io}^\star =\$0,
\forall i\in I.  
\end{aligned}  \end{equation}
 \vspace{-1em}
\begin{equation} \begin{aligned}
\relax [\sum_{\forall j\in \mathcal{E}_o^{VGA}} (y_{rj}\pi_{jo}^\star )- y_{ro}(1+P_{ro}^\star )]u_{ro}^\star =\$0,
 \forall r\in R .
 \label {Eq:36} \end{aligned} \end{equation}
\indent $v_o$ and $u_o$ are free-in-signs in the TVG programs. Since \((\textit{X, Y})>0\) and \(\tau_o^\star >\$0\) yields the estimated virtual unit prices, \(v_o^\star  >0, u_o^\star >0\). Therefore, the braces at the left-hand side of  (\ref{Eq:35}) and (\ref{Eq:36}) equal zero. 

\indent For the PT model, each DMU-j has the property as shown in  (\ref{Eq:37}). When $(v_o^\star x_j - u_o^\star y_j )=0$, DMU-j belongs to \(\mathcal{E}_o^{PT}\)  and \(\pi_{jo}^\star >0\). The \emph{Efficiency Equator} is the diagonal line at the origin in Figure \ref{fig:2}, that \((v_o^\star x_j - u_o^\star y_j)=\$0, \forall j\in \mathcal{E}_o^{PT}\).
\begin{equation}
	(v_o^\star x_j - u_o^\star y_j )\pi_{jo}^\star = \$0, \forall j\in J .
	\label{Eq:37}
\end{equation} 
For the TSc model, each DMU-j has the property shown in  (\ref{Eq:38}). $w_o^\star$  is the estimated unit price of the SIC scalar, and $\kappa_o^c w_o^\star$  is the vScalar. When $(v_o^\star x_j - u_o^\star y_j + 1 w_o^\star ) = \$0$, DMU-j belongs to $\mathcal{E}_o^{TS}$ and $\pi_{jo}^\star>0$. Otherwise, $(v_o^\star x_j - u_o^\star y_j + 1  w_o^\star ) > \$0$, DMU-j does not belong to $\mathcal{E}_o^{TS}$ and $\pi_{jo}^\star=0$. The Efficiency Equator is the diagonal line at the origin in Figure \ref{fig:3}, that $(v_o^\star x_j - u_o^\star y_j + 1  w_o^\star)=\$0, \forall j\in \mathcal{E}_o^{TS}$. 
 \vspace{-.5em}
\begin{equation}
	 (v_o^\star x_j - u_o^\star y_j + w_o^\star ) \pi_{jo}^\star = \$0, \forall j\in J.
	 	\label{Eq:38} \end{equation}

 (\ref{Eq:39}) depicts the relationships between the estimated decision variables in the TVG and TAP programs.
\begin{equation}\begin{aligned}
	(v_{io}^\star x_{io} - \tau_o^\star) Q_{io}^\star = \$0,\forall i \in I ; (u_{ro}^\star y_{ro} - \tau_o^\star) P_{ro}^\star= \$0,\forall r\in R .\end{aligned} \label{Eq:39} \end{equation}
 In  (\ref{Eq:40}), $\sum_{\forall j\in \mathcal{E}_o^{TS}}\pi_{jo}^\star - \kappa_o^c = 0$ so that TVG of the TSc model is decreased, constant, and increased when $w_o^\star>\$0,  w_o^\star=\$0, \textrm{ and } w_o^\star<\$0$, respectively. 
\begin{equation} \label{Eq:40}
	(\sum_{\forall j\in \mathcal{E}_o^{TS}} \pi_{jo}^\star -\kappa_o^c ) w_o^\star  = \$0.			\end{equation}
\subsection{Post-Analysis of the Optimal Solutions}\label{sec:2.5}
The optimal solutions of the PT and TSc VGA models can quantify the assessment items. 
\subsubsection{Virtual Technology Sets}\label{sec:2.5.1}
\indent The formulation and definitions of the \textit{virtual technology set} in the context of Step II within the PT and TSc models are presented in  (\ref{Eq:41}) and (\ref{Eq:42}); Step I has the similar expressions. In which, the virtual scales of each DMU-j (\textit{pvInput}, \textit{pvOuput}) and (\textit{avInput}, \textit{avOutput}) are defined in (\ref{Eq:12}) and (\ref{Eq:28}). The virtual scales in Step I for each DMU-j represent certain aspects or characteristics within the model. These sets likely encapsulate specific parameters, constraints, or variables pertinent to the subsequent steps' formulation and computation. Virtual scales in Step II might represent refined or adjusted versions of the initial sets of Step I, potentially incorporating the outcomes or adjustments derived from the earlier steps.

\indent The transition from Step I to Step II usually involves refining or re-calibrating the parameters or sets based on the intermediate solutions or computations from the preceding step. This iterative process often helps converge toward more accurate or optimal results within the models.
\begin{equation}
\begin{aligned}
	\Phi_o^{PT\star}&=\lbrace (\alpha^\star, \beta^\star) \mid (\alpha_j^{\star},\beta_j^{\star}), \forall j\in J\rbrace. 
	\end{aligned} 
	\label{Eq:41}
	\end{equation}
 \vspace{-1.5em}
\begin{equation} \begin{aligned} \label{Eq:42}  
 \Phi_o^{TSc\star}=
\lbrace (\alpha^{c\star},\beta^{c\star})\mid(\alpha_j^{c\star}, \beta_j^{c\star} ), \forall j\in J& \rbrace. 
	 \end{aligned}
	 \end{equation}	 
\subsubsection{Return-to-virtual-scale (RTvS)}\label{sec:2.5.2} ($\ref{Eq:43}$) and ($\ref{Eq:44}$) depict the computation steps for DMU-o's \textit{benchmark virtual scales} and \textit{affected benchmark virtual scale} in the PT and TSc models during Step II,  ( bvInput, bvOutput) =  ($\widehat\alpha_o^\star$, $\widehat\beta_o^\star$) and ( abvInput,  abvOutput) =  ($\widehat\alpha_o^{c\star}$, $\widehat\beta_o^{c\star}$), respectively. These equations likely involve transforming or deriving benchmark values based on specific conditions or constraints within the models.

\indent The expressions for \textit{bvInput, bvOutput, abvInput, and abvOutput} in Step I are similar to those in  (\ref{Eq:43}) and (\ref{Eq:44}) with the replacement of the superscript "$\star$" with "$\#$." This step-wise process of deriving and refining benchmarks could indicate a progressive approach where the initial estimates from Step I are further modified or calibrated in Step II based on additional considerations or parameters introduced in the models.
\begin{equation} \label{Eq:43} 
	\widehat{\alpha}_o^\star= \sum_{\forall i\in I} \widehat{x}_{io} v_{io}^\star \textrm{  and  } \widehat{\beta}_o^\star= \sum_{\forall r\in R}\widehat{y}_{ro} u_{ro}^\star.	
	\end{equation}	
  \vspace{-1.5em}
\begin{equation}\label{Eq:44} \begin{aligned}	
	\widehat\alpha_o^{c\star} = [v_o^\star \widehat {x}_o + (1 - \gamma_o)\omega_o^{c\star} ]; \widehat\beta_o^{c\star} = (u_o^\star  \widehat{y}_o - \gamma_o\omega_o^{c\star} ).
\end{aligned}
\end{equation}
\indent  ($\ref{Eq:45}$) and ($\ref{Eq:46}$) compute  $\Xi_o^{PT\star}$  and  $\Xi_o^{TSc\star}$ that based on the virtual scales. The RTvS values might represent a measure or index that assesses the efficiency or productivity performance of DMU-o within the context of these specific models.
\begin{equation}
	\Xi_o^{PT\star}  = (\widehat\beta_o^\star/\beta_o^\star ) / (\widehat\alpha_o^\star/\alpha_o^\star)=1/E_o^{PT\star}.
	\label{Eq:45}
	\end{equation}
  \vspace{-1em}
\begin{equation} 
 \label{Eq:46}
	\Xi_o^{TSc\star} = (\widehat\beta_o^{c\star} /\beta_o^{c\star}) / (\widehat\alpha_o^{c\star} /\alpha_o^{c\star} )
  =1/E_o^{TSc\star}.
\end{equation}
\indent For the case that $\omega_o^\star\geq \$0$, $\Xi_o^{TSc\star}$ and $E_o^{TSc\star}$ will be decreased and increased as the SIC scalar $\kappa_o^c$ rises. Conversely, when $\omega_o^\star\leq \$0$, $\Xi_o^{TSc\star}$ and $E_o^{TSc\star}$ will be decreased and increased, respectively, as the SIC scalar $\kappa_o^c$ diminishes. A preference is given to higher $\Xi_o^{TSc\star}$ and $E_o^{TSc\star}$. Depending on whether $\omega_o^\star\geq \$0$ or $\omega_o^\star\leq \$0$, DMU-o has the flexibility to adjust $\kappa_o^c$ towards its lower or upper bound. This adjustment aims to achieve the final $\kappa_o^z$  that encompasses achievable benchmarks for inputs and outputs. The implications of these dynamics are further elucidated through the numerical examples and graphical illustrations presented in Section \ref{sec:5.1}, which depict how varying $\kappa_o^c$  influences the efficiency and the position of DMU-o relative to the Efficiency Equator.
\subsubsection{Interconnections Between Inputs and Outputs Indices} \label{sec:2.5.3}  ($\ref{Eq:47}$) and ($\ref{Eq:48}$) demonstrate adjustments to the affected virtual prices of an input, $(v_{io}^\star x_{io}  + \gamma_{io}^Q \omega_o^{c\star} )$, and an output, $(u_{ro}^\star y_{ro} - \gamma_{ro}^P \omega_o^{c\star} )$, respectively considering the proportions to the vScalar, represented by  $\gamma_{io}^Q$  and  $\gamma_{ro}^P$. 
\begin{equation}\begin{aligned}
	\Delta_o^{TSc\star} 
	= \sum_{\forall i\in I} (v_{io}^\star x_{io} + \gamma_{io}^Q \omega_o^{c\star} )  - \sum_{\forall r\in R}(u_{ro}^\star y_{ro} - \gamma_{ro}^P \omega_o^{c\star} ).
\end{aligned}	\label{Eq:47}
	\end{equation}	
 \vspace{-1.5em}
\begin{equation} \begin{aligned}
	\gamma_{io} ^Q  &= (1 - \gamma_o)Q_{io}^\star/\sum_{\forall i\in I} Q_{io}^\star, \forall i \in I ;\\
 \gamma_{ro}^P &= \gamma_o P_{ro} ^\star / \sum_{\forall r\in R} P_{ro}^\star, \forall r \in R.
	\end{aligned}\label{Eq:48}
	\end{equation}
\indent In  ($\ref{Eq:47}$) and ($\ref{Eq:48}$), the proportions involving ($v_{io}^\star x_{io} +  \gamma_{io}^Q \omega_o^{c\star}$ ) and ($u_{ro}^\star y_{ro}-\gamma_{ro}^P \omega_o^{c\star} $) indeed highlight the interconnectedness between the input and output indices. These proportions indicate the adjustments made to the affected virtual prices concerning the quantities represented by $\gamma_{io}^Q$  and $\gamma_{ro}^P$, displaying how changes in one influence the other within the context of the model variables and constraints.
\subsubsection {2D Graphic Intuitions}\label{sec:2.5.4}
DMU-o can visualize the outcomes of the post-analysis based on the optimal solutions from the VGA models. For instance, I am using the example dataset in Table \ref{table:1}. One can identify the third SIC scalar $\kappa_o^3$ by try and error. While the  $\kappa_o^3$ may not be within the range of  $\kappa_o^1$ and  $\kappa_o^2$. The particular $\kappa_o^3$ within the TS3 model produces $E_o^{PT\star}=E_o^{TS3\star}$. Pure technical efficiency is affected by the SIC. 
Figure \ref{fig:2} depicts the solutions of PT and TS3 models. Points Kp, K3, AP3, and O are located at ($\alpha_o^{\star},\beta_o^{\star}$), ($\alpha_o^{3\star}, \beta_o^{3\star} $), ($[1-\gamma_o] \omega_ o^{3\star }, -\gamma_o  \omega_o^{3\star}$), and (0,0), respectively. Points O, K3, and AP3 at the triangle. The upper-right of the figure shows the two rectilinear distances between the pair points (Kp, Tp) and (K3, T3), the virtual gaps estimated by the PT and TS3 models. DMU-o could comprehend the results between the PT and TS3 models to understand the effects of the SIC.

\indent Similarly, in Figure \ref{fig:3}, Points O, K1, and AP1 of the triangle are solutions of the TS1 model. The three points O, K2, and AP2 at the triangle are solutions of the TS2 model. Figure \ref{fig:3} showcases the effects of using different SIC scalars within a certain interval and helping select a final scalar $\kappa_o^z$ for productivity management. In this example, due to $w_o^\star >0$, we have $E_o^{TS2\star}$ $>$ $E_o^{TSz\star}$ $>$ $E_o^{TS1\star}$. The TSz model provides the most preferred AARs ($Q_o^{TSz\star}, P_o ^{TSz\star}$).

In Figure \ref{fig:3}, Points O, APc, and Kc at the triangle which expresses solutions of the TSc model, where "c" denotes the choices 1, 2, and 3. (\ref{Eq:49}), (\ref{Eq:50}), and (\ref{Eq:51}) are the expressions of the slopes for the vectors $\overrightarrow{\rm O, Kc}$, $\overrightarrow{\rm O, APc}$, and $\overrightarrow{\rm APc, Kc}$, respectively. The ratio of the scalar prices of inputs to outputs equals the slope $\bar{m}(\overrightarrow{\rm O, APc})$. Similarly, the ratio of the gray virtual input to output equals the slope $\bar{m}(\overrightarrow{\rm APc, Kc})$.  
 \begin{equation} \begin{aligned}
	\bar{m}(\overrightarrow{\rm O,Kc})= \frac{{\beta_o^{c\star}-0}}{{\alpha_o^{c\star}-0}} =\frac{{\beta_o^{c\star}}}{{\alpha_o^{c\star}}}.
	\label{Eq:49}
	\end{aligned} \end{equation}
 \vspace{-1em}
\begin{equation}\begin{aligned}
	 \bar{m}(\overrightarrow{\rm O,APc})=\frac{-\gamma_o \omega_ o^{c\star }-0}{(1-\gamma_o) \omega_ o^{c\star }-0}=\frac{-\gamma_o \omega_ o^{c\star }}{(1-\gamma_o) \omega_ o^{c\star }}.
	\end{aligned}\label{Eq:50}
	\end{equation}	
 \vspace{-1em}
 \begin{equation}\begin{aligned}
	\bar{m}(\overrightarrow{\rm APc,Kc})=\frac{\beta_o^{c\star}-(-\gamma_o \omega_ o^{c\star })}{{\alpha_o^{c\star}}-[-(1-\gamma_o) \omega_ o^{c\star }]}
 =\frac{u_o^\star y_o}{ v_o^\star x_o}.
	\end{aligned}\label{Eq:51}
	\end{equation}	
 For instance, in Figure \ref{fig:2}, the triangle of the TS3 model has the following expressions. 
 \begin{equation} 
\begin{aligned}
	& \overrightarrow{\rm O,K3} =\overrightarrow{\rm O,AP3}+\overrightarrow{\rm AP3,K3} 
	 \iff \frac{\beta_o^{3\star}}{\alpha_o^{3\star}} =\frac{ u_o^\star y_o -\gamma_o  \omega_o^{3\star}}{ v_o^\star x_o +(1-\gamma_o)\omega_o^{3\star} }.
 \end{aligned} \label{Eq:52}
	 \end{equation}
\begin{theorem}\label{thm:5}{ TSc model fulfills Condition C3.}
\end{theorem}

\noindent\textbf{Proof.} As shown in (\ref{Eq:53}), the estimated TSc efficiency score has the following expressions. 
\begin{equation} \begin{aligned}
 \overrightarrow{\rm O,Kc} =\overrightarrow{\rm O,APc}+\overrightarrow{\rm APc,Kc};\quad 
 E_o^{TSc\star}=\bar{m}(\overrightarrow{\rm O,Kc})= \frac{\beta_o^{c\star}}{\alpha_o^{c\star}}.
\label{Eq:53} \end{aligned}
\end{equation}
The effects of the SIC within the TSc model could be visualized as triangular with points O, APc, and Kc in the 2D graphical intuition. The slopes $\bar{m}(\overrightarrow{\rm O,APc})$ and $\bar{m}(\overrightarrow{\rm APc, Kc})$ are not the scale and technical efficiencies. \hfill $\square$

\section{Second Scenario: Evaluating Efficient DMUs' Super-Efficiencies}\label{sec:3}
In the first scenario, DMU-o is identified as efficient when $E_o^{sPT\star}$=1. Let efficient DMUs belong to the set \textbf{\.{E}}. PT and TSc models are modified into sPT and sTEc models to estimate \textit{all allowance ratios (AARs)} of each DMU belonging to \textbf{\.{E}}, named DMU-o. DMU-o has a super-efficiency larger than one and will deteriorate to 1 when DMU-o expands inputs and deduces outputs. A DMU is superior to other efficient DMUs as it has a higher super-efficiency. Variables $Q_{io}$ and $P_{ro}$ are changed to denote the allowance ratios for DMU-o. In particular, DMU-o is excluded from the reference set $J$ when assessing its super-efficiency. 
\subsection{Pure Technical Super-Efficiency (sPT) Model}\label{sec:3.1}
Total adjustment price (TAP), the dual program of the sPT model:	
\begin{equation}
	\delta_o^{sPT\star} = \min_{Q_o, P_o, \pi_o} \sum_{\forall i\in I} Q_{io} \tau_o  + \sum_{\forall r\in R}P_{ro} \tau_o, \forall o \in \textbf{\.{E}};
	\label{Eq:54}
 \end{equation} \vspace{-1em}
\begin{equation}
	s.t. -\sum_{\forall j\in J-\{o\}} x_{ij} \pi_{jo} + Q_{io} x_{io} \geq - x_{io}  , \forall i\in I;
 \label{Eq:55}
 \end{equation}  \vspace{-1.5em}
\begin{equation} 
 \qquad \sum_{\forall j \in J-\{o\}} y_{rj} \pi_{jo} + P_{ro} y_{ro} =  y_{ro}, \forall r \in R;	
\label{Eq:56}
 \end{equation}  \vspace{-2em}
\begin{equation} \label{Eq:57}
\qquad \pi_o, Q_o  , P_o  \ge 0.	
 \end{equation}
\noindent Total virtual gap (TVG), the primal program of the sPT model:
\begin{equation} \label{Eq:58}
	\Delta_o^{sPT\star} = \max_{v_o, u_o} -\sum_{\forall i\in I} v_{io} x_{io}  + \sum_{\forall r\in R} u_{ro} y_{ro}, \forall o \in \textbf{\.{E}};
 \end{equation}  \vspace{-1em}
 \begin{equation} \begin{aligned} 
 \label{Eq:59}
	s.t. - \sum_{\forall i\in I} v_{io} x_{ij}  &+ \sum_{\forall r\in R}u_{ro}  y_{rj}  \le 0, \forall j\in J-\{o\}; 
 \end{aligned}
 \end{equation}  \vspace{-2em}
 \begin{equation}
	x_{io} v_{io} \ge\tau_o, \forall i\in I;	
 \label{Eq:60}
 \end{equation}  \vspace{-2 em}
\begin{equation}
	y_{ro} u_{ro}\ge\tau_o,\forall r\in R;
\label{Eq:61}
 \end{equation} \vspace{-2 em}
\begin{equation}
	v_o\geq 0\quad and \quad u_o \quad free.	
 \label{Eq:62}
 \end{equation}
 (\ref{Eq:59}) limits each DMU-j, except DMU-o, upper bound to zero virtual gaps. While (\ref{Eq:60}) and (\ref{Eq:61}) restrict virtual prices of DMU-o lower bound to the unified goal price $\tau_o$. 

\indent Use Step I and Step II of the sPT model to determine the goal price, as shown in (\ref{Eq:63}). The estimated $u_{o}^\#  y_o$ equals \(\$\)1 and larger than $v_{o}^\#  x_o$.
\begin{equation}
	\bar{t} = \$ 1/u_{o}^\#  y_o \textrm {  and  } \tau_o^\star   = \$ \bar{t}. 
	\label{Eq:63}
  \end{equation}
Similar to (\ref{Eq:11}) and (\ref{Eq:12}), the optimal solutions of  (\ref{Eq:54}) and (\ref{Eq:58}) as depicted in (\ref{Eq:64}) and (\ref{Eq:65}).
The \textit{total virtual gaps} (TVG) are equal to  $(-\alpha_o^\# +\beta_o^\#)$ and  $(-\alpha_o^\star +\beta_o^\star  )$ in Steps I and-II .  
\begin{equation}  \label{Eq:64} 
 \begin{aligned}
& \$0 <\delta_o^{sPT\#} = \sum_{\forall i\in I} Q_{io} ^\# \times \$1  + \sum_{\forall r\in R}P_{ro} ^\#  \times \$1 =\Delta_o^{sPT\#} ; \\ 
&\$ 0 <\delta_o^{sPT\star } = \sum_{\forall i\in I} Q_{io}^\star  \tau_o ^\star  + \sum_{\forall r\in R} P_{ro}^\star  \tau_o ^\star=\delta_{xo}^{sPT\star} + \delta_{yo}^{sPT\star} =\Delta_o^{sPT\star}< \$1.
  \end{aligned} 
 \end{equation} 
 \begin{equation}
\begin{aligned}
&\$0 < \Delta_j^{sPT\#}  =- v_o^\# x_j + u_o^\# y_j = -\alpha_j^\# + \beta_j^\# ;\\
&\$0 <\Delta_j^{sPT\star }  = -v_o^\star  x_j + u_o^\star  y_j =-\alpha_j^\star  + \beta_j^\star ;
j \in J -\{o\}.
\end{aligned}
\label {Eq:65}
 \end{equation}
\begin{theorem}\label{thm:6}
{The sPT model estimated sPT efficiency no less than 1.}
\end{theorem}
\noindent\textbf{Proof.} The sPT model estimates the technical parameters $v_o^\star$ and $u_o^\star$ to have the maximum virtual gap  $\Delta_o^{sPT\star }$. Which is converted into the $\textit{super-pure-technical-efficiency}$ (sPT)    $E_o^{sPT\star }$ according to (\ref{Eq:66}). 
\begin{equation}
\begin{aligned} \label{Eq:66} 
	&1\leq E_o^{sPT\star } = \beta_o^\star  /\alpha_o^\star.	  
\end{aligned}
\end{equation}
$E_o^{sPT\star }$ is no less than one,  where $\alpha_o^\star \leq$ $\beta_o^\star = \$1$. \hfill $\square$\\
 \indent (\ref{Eq:67}) likely represents this normalization process, where the solutions obtained in Step I are adjusted or transformed to achieve the solutions in Step II. 
\begin{equation} \begin{aligned}
  (\delta_o^{sPT\star },  \Delta_o^{sPT\star },  v_o^\star, u_o^\star)
  = \bar{t} \times (\delta_o^{sPT\#}, \Delta_o^{sPT\#}, v_o^\#, u_o^\#).
\label {Eq:67} \end{aligned}
\end{equation}
The symbol  $\mathcal{E}_o^{sPT}$  denotes the reference set of DMU-o in the sPT evaluations. Each peer DMU-j belonging to the reference set has sPT efficiency equals 1, and its estimated intensity, $\pi_{jo}^\star >0$. At the same time, the estimated intensity of the other inefficient DMUs, $\pi_{jo}^\star =0$. DMU-o estimates the benchmark of each input-i and output-r,  $\widehat{x}_{io}^{sPT\star}$ and  $\widehat{y}_{ro}^{sPT\star}$ via  (\ref{Eq:68}). DMU-o mimics the reference DMUs with their estimated intensities,  $\pi_{jo}^\star$, $\forall j\in \mathcal{E}_o^{sPT}$. DMU-o is superior to its reference peers, whose criteria configurations are analog to DMU-o.

\begin{equation}
\begin{aligned}
	&\widehat{x}_{io}^{sPT\star} = \sum_{\forall j\in \mathcal{E}_o^{sPT}} x_{ij}  \pi_{jo}^\star  = x_{io} (1 + Q_{io}^\star ), \forall i\in I; \\
	&\widehat{y}_{ro} ^{sPT\star} = \sum_{\forall j\in \mathcal{E}_o^{sPT}}y_{rj} \pi_{jo}^\star =y_{ro} (1 - P_{ro}^\star  )
 \forall r\in R.
\end{aligned}
\label {Eq:68}
\end{equation}
Assessing DMU-o with the adjusted input-i and output-r, \(\widehat{x}_{io}^{sPT\star}\) and \(\widehat{y}_{ro} ^{sPT\star} \), as shown in (\ref{Eq:68}), will have the sPT efficiency   $\widehat{E}_o^{sPT\star }$ equals 1. 
\begin{equation}
    \sum_{\forall j\in \mathcal{E}_o^{sPT}}\pi_{jo}^\star   = \kappa_o^1 . 
 \label {Eq:69} \end{equation}
Let the total of estimated intensities be the first SIC scalar, $\kappa_o^1$, computed via  (\ref{Eq:69}). 

\subsection{Super-Efficiency sTSc Model}\label{sec:3.2}
 Adding the SIC (\ref{Eq:73}) to the sPT model's dual program will have the next program. The SIC corresponds to the free-in-sign decision variable, $w_o^c$.
\noindent TAP (dual) program of the sTSc model:	
\begin{equation}\label {Eq:70} 
	\delta_o^{sTSc\star}= \min_{Q_o, P_o} \sum_{\forall i\in I} Q_{io} \tau_o + \sum_{\forall r\in R}P_{ro}   \tau_o, \forall o \in \textbf{\.{E}};	
\end{equation}  \vspace{- 1.5 em}
 \begin{equation} \label{Eq:71}
	s.t. -\sum_{\forall j\in J-\{o\}} x_{ij} \pi_{jo} + Q_{io} x_{io}  \geq - x_{io}, \forall i \in I;	
 \end{equation}  \vspace{- 1.5 em}
\begin{equation}\label {Eq:72}
	 \sum_{\forall j\in J-\{o\}} y_{rj} \pi_{jo} + P_{ro} y_{ro} =  y_{ro}  , \forall r \in R;	
\end{equation}  \vspace{- 1.5 em}
 \begin{equation}\label {Eq:73}
	 \sum_{\forall j\in J-\{o\} } \pi_{jo}=\kappa_o^c;
  \end{equation}  \vspace{- 1.5 em}
  \begin{equation}\label {Eq:74}
	\pi_o,  Q_o, P_o  \ge 0.	
\end{equation}
  TVG (primal) program of the sTSc model:
\begin{equation}\begin{aligned}
	\Delta_o^{sTSc\star} =  \max_{x_o, v_o, w_o^c }  -\sum_{\forall i\in I} v_{io} x_{io}  &+ \sum_{\forall r\in R} u_{ro} y_{ro}  + \kappa_o^c w_o^c, \forall o \in \textbf{\.{E}};
 \label {Eq:75}
 \end{aligned}
 \end{equation}
  \vspace{- 1.5 em}
 \begin{equation} \begin{aligned}\label {Eq:76}
s.t. -\sum_{\forall i\in I} v_{io} x_{ij} + \sum_{\forall r\in R} u_{ro} y_{rj}  + 1 w_o^c  \le 0,
\forall j \in J-\{o\}; 
\end{aligned} \end{equation}  \vspace{- 2 em}
 \begin{equation}\label{Eq:77}
x_{io} v_{io} \ge\tau_o,\forall i \in I; 
\end{equation}  \vspace{-2em}
 \begin{equation}\label {Eq:78}
y_{ro} u_{ro} \ge \tau_o, \forall r \in R;	
 \end{equation}  \vspace{-2em}
 \begin{equation}\label {Eq:79}
v_o\geq 0, u_o, w_o^c \quad free.
\end{equation}

\indent Similar to the TSc model, the sTSc model determines the goal price in two steps.  \(\tau_o^\# = \$ 1\) and \(\tau_o^\star   = \$ \bar{t}\) will solve the TAP program comprehensively with the solutions \((Q_o^\star , P_o^\star , \pi^\star )=(Q_o^\#, P_o^\#, \pi^\# )\). Use (\ref{Eq:27}) to obtain \(\gamma_o\) and \((1- \gamma_o)\). Partition the \textit{vScalar} into two parts to reflect the effects of the SIC on inputs and outputs of DMU-o.

\indent The maximized TVG of DMU-o in Step II of (\ref{Eq:75}) is expressed as the \textit{affected virtual Output (avOutput)} minus \textit{affected virtual Input (avInput)}, as shown in (\ref{Eq:80}). The solutions in Step I have a similar expression.
\begin{equation}\begin{aligned} \label{Eq:80}
&\$0 < \Delta_o^{sTSc\star} = -v_o^\star x_o+u_o^\star y_o+\omega_o^{c\star }\\
&= -[v_o^\star x_o-(1-\gamma_o)\omega_o^{c\star }]+(u_o^\star   y_o+\gamma_o\omega_o^{c\star }) \\
& =-[gvInput-ivScalar]+[gvOutput+ovScalar]\\
	 &= -avInput+avOutput= -\alpha_o^{c\star} + \beta_o^{c\star}< \$ 1;\\
& \$0 <\Delta_o^{sTSc\#} =-\alpha_o^{c\#} + \beta_o^{c\#}.
 \end{aligned}  \end{equation}
 \indent Similarly, the estimated \textit{vScalar} of Step II  in  (\ref{Eq:76}), \(1w_o^*\), is decomposed into two components, as shown in  (\ref{Eq:81}); the solutions of Step I have the similar expressions. 
  \begin{equation}\begin{aligned} \label{Eq:81}
&\$0 < \Delta_j^{sTSc\star}=- v_o^\star x_j+u_o^\star y_j+ 1w_o^{c\star } \\
&= -[v_o^\star x_j-(1-\gamma_o) 1w_o^{c\star }]+(u_o^\star   y_j+\gamma_o 1w_o^{c\star }) = -\alpha_j^{c\star } + \beta_j^{c\star }, \forall j \in J,  j\neq {o} .
 \end{aligned} \end{equation}
 Each DMU-j is expressed by the pair of \emph{(avInput, avOutput)}, the \textit{virtual scales}. The symbol  \(\mathcal{E}_o^{sTS}\) denotes the set of reference peers in evaluating DMU-o by the sTS model, where each reference peer DMU-j has \(\pi_{jo}^\star >0\). The other inefficient DMU-j has \(\pi_{jo}^\star =0\). In  (\ref{Eq:81}), any best DMU-j belongs to \(\mathcal{E}_o^{sTS}\) has $\alpha_j^{c\#} = \beta_j^{c\#}$ and  $\alpha_j^{c\star } = \beta_j^{c\star }$, while the other remaining DMUs have $\alpha_j^{c\# } > \beta_j^{c\# }$ and $\alpha_j^{c\star } > \beta_j^{c\star }$. 
\begin{theorem}\label{thm:7} {The sTSc model estimates sTSc efficiency score larger than 1.}
\end{theorem}
\noindent\textbf{Proof.} Using (\ref{Eq:82}) would have the dimensionless value of $\bar{t}$. 
\begin{equation} \begin{aligned} \label {Eq:82}
	\$\bar{t }= \$ 1 / \beta_o^{c\#}= \$ 1 / (u_o^\#  y_o + \gamma_o \omega_o^{c\#})
 =\tau_o^\star.
\end{aligned} \end{equation}
Because $\tau_o^\#:\tau_o^\star  = \$1 : \$\bar{t}$, therefore, if $\tau_o^\# $ = $\$\bar{t}$, then \(\beta_o^{c\star}\) equals \$1. According to (\ref{Eq:80}), the avInput and avOutput, $\alpha_o^{c\star}$ and $\beta_o^{c\star}$, contain the gvInput, gvOutput, ivScalar, and ovScalar. $-\alpha_o^{c\star}$ +  $\beta_o^{c\star}$ is larger than \$0. Using (\ref{Eq:83}), the estimated virtual gap  $\Delta_o^{sTSc\star }$ is converted into the maximized sTSc super-efficiency score,  $E_o^{sTSc\star}$, should be larger than 1. 
\begin{equation}
\begin{aligned}
	E_o^{sTSc\star } =\frac{
	 u_o^\star y_o+\gamma_o \omega_o^{c\star}}{v_{o}^\star x_o-(1-\gamma_o)\omega_o^{c\star} } = \frac{\beta_o^{c\star}}{\alpha_o^{c\star}} > 1.	\end{aligned}
 \label {Eq:83}
  \end{equation}
 \hfill $\square$ 

\indent Normalizing the solutions of Step I by the dimensionless value $\bar{t}$ would obtain the solutions of Step II, as shown in  (\ref{Eq:84}). 
\begin{equation}\begin{aligned}\label {Eq:84}
	(\Delta_o^{sTSc\star }, \delta_o^{sTSc\star }, v_o^\star , u_o^\star , w_o^{c\star })
	=  \quad\bar{t}\times(\Delta_o^{sTSc\#},\delta_o^{sTSc\#} ,v_o^\#, u_o^\# , w_o^{c\#}).
\end{aligned}
\end{equation}
\indent Use  (\ref{Eq:85}) to compute the benchmark of each performance index,  \(\widehat{x}_{io}^{sTSc\star}\) , and  \(\widehat{y}_{ro}^{sTSc\star}\). DMU-o imitates the reference peers with their estimated intensities \((\pi_{jo}^\star )\). 
\begin{equation}
\begin{aligned}\label {Eq:85} 
	\widehat{x}_{io}^{sTSc\star} &= \sum_{\forall j\in \mathcal{E}_o^{sTS}}x_{ij}\pi_{jo}^\star = x_{io}(1 + Q_{io}^\star ), \forall i\in I;\\
	\widehat{y}_{ro} ^{sTSc\star} &= \sum_{\forall j\in \mathcal{E}_o^{sTS}}y_{rj}\pi_{jo}^\star = y_{ro}(1 - P_{ro} ^\star ), \forall r\in R .   \end{aligned} \end{equation}

\begin{theorem}\label{thm:8}{DMU-o has an sTSc super-efficiency equals 1 when using the estimated benchmarks.}
\end{theorem}
\noindent\textbf{Proof.} 
DMU-o uses the benchmarks, $\widehat{x}_{io}^{TSc\star}$ and $\widehat{y}_{ro} ^{TSc\star}$, will decrease the sTSc super-efficiency $\widehat{E}_o^{sTSc\star }$ to 1 while the $\delta_o^{TSc\star}$ drops to $\$$0. \hfill $\square$
\subsection{The SIC Effects in Measuring the Super-Efficiency}\label{sec:3.3}
In Figure \ref{fig:4}, the coordinates of points O, AP1, AP2, B1, and B2 of the sTS1 and sTS2 models in evaluating DMU-B can be visualized. (\ref{Eq:86}),(\ref{Eq:87}), and (\ref{Eq:88}) compute the slopes of the vectors $\overrightarrow{\rm O, Bc}$, $\overrightarrow{\rm O, APc}$, and $\overrightarrow{\rm APc, Bc}$, respectively, where the choice "c" could be "1" and "2."  
\begin{equation} \begin{aligned}
	\bar{m}(\overrightarrow{\rm O,Bc})= \frac{{\beta_o^{c\star}-0}}{{\alpha_o^{c\star}-0}} =\frac{{\beta_o^{c\star}}}{{\alpha_o^{c\star}}}.
	\label{Eq:86}
	\end{aligned} \end{equation}	
 \vspace{-1.5em}
\begin{equation}\begin{aligned}
	 \bar{m}(\overrightarrow{\rm O,APc})=\frac{\gamma_o \omega_ o^{c\star }-0}{-(1-\gamma_o) \omega_ o^{c\star }-0}=\frac{\gamma_o \omega_ o^{c\star }}{-(1-\gamma_o) \omega_ o^{c\star }}.
	\end{aligned}\label{Eq:87}
	\end{equation}	
 \vspace{-1em}
 \begin{equation}\begin{aligned}
	\bar{m}(\overrightarrow{\rm APc,Bc})=\frac{\beta_o^{c\star}-\gamma_o \omega_ o^{c\star }}{{\alpha_o^{c\star}}-[-(1-\gamma_o) \omega_ o^{c\star }]}=\frac{u_o^\star y_o}{v_o^\star x_o}.
	\end{aligned}\label{Eq:88}
	\end{equation}	
 The same analysis should apply to DMU-D, as Figure \ref{fig:4} depicts.
\begin{theorem}\label{thm:9}{ The sTSc model fulfills Condition C3 that the relative efficiency are affected by the SIC.}
\end{theorem}

\noindent\textbf{Proof.} The three vectors have the relationships as shown in (\ref{Eq:89}). The estimated sTSc super-efficiency score is affected by the SIC. 
\begin{equation} \begin{aligned}
 \overrightarrow{\rm O,Bc} =\overrightarrow{\rm O,APc}+\overrightarrow{\rm APc,Bc};\quad
 E_o^{sTSc\star}=\bar{m}(\overrightarrow{\rm O,Bc})= \frac{\beta_o^{c\star}}{\alpha_o^{c\star}}.
\label{Eq:89} \end{aligned}
\end{equation}
The effects of the SIC within the sTSc model could be visualized as triangular with points O, APc, and Bc in the 2D graphical intuition. The slopes $\bar{m}(\overrightarrow{\rm O,APc})$ and $\bar{m}(\overrightarrow{\rm APc, Bc})$ are not the scale and technical efficiencies. \hfill $\square$

\section{Four-phase VGA Procedure}\label{sec:4}
The following presentation is to calibrate the AARs to measure the inefficiency of DMU-o. It can be applied to measure super-efficiency. Nonetheless, the PT model estimated AARs may not be achievable for DMU-o. The TSc model with the SIC equates to the choice scalar, $\kappa_o^c,$ that identifies the reference peers in $\mathcal{E}_o^{TS}$. The benchmarks of inputs and outputs are the function of intensities of the reference peers belonging to $\mathcal{E}_o^{TS}$, as shown in (\ref{Eq:33}). Section \ref{sec:2.5.4} illustrates the effects of SIC on the solutions.    

\begin{figure}[H] \centering
 \includegraphics
 [width=0.6\linewidth, height=8cm]{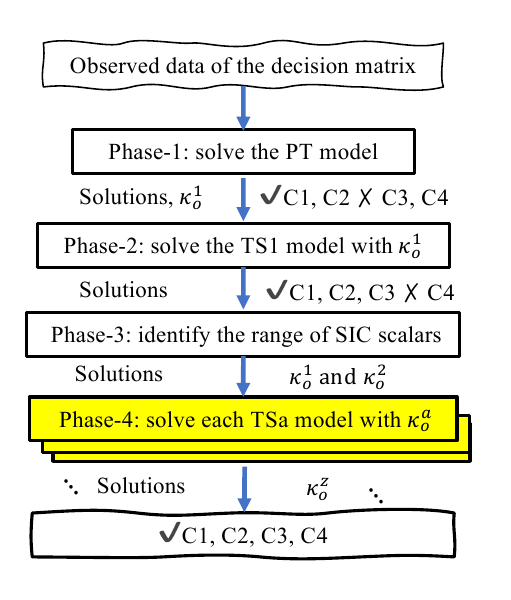} 
 \caption {Four-phase VGA procedure. } \label{fig:1}
  \end{figure}

\indent The four-phase process illustrated in Figure \ref{fig:1} significantly improves the productivity management of DMUs with Conditions C1, C2, C3, and C4.

In summary, the four conditions collectively guide the estimation of AARs by recognizing efficiency's inherent and operational dimensions, promoting continuous improvement, and setting realistic targets for DMU-o. 

This framework is structured around four conditions that epitomize a performance improvement problem. It incorporates practical considerations into the PT and TSc VGA models for systematic solutions. By integrating innovations that cater to the third and fourth conditions, our approach offers a more comprehensive and systematic method for efficiency assessment with potential applicability across various sectors.
\subsection{Phase 1: PT model Identifies the First Scalar.} \label{sec:4.1} 
An initial assessment distinguishes between efficient and inefficient DMUs in the performance evaluation. This distinction is crucial for the subsequent analysis and is primarily based on the existence of a "virtual gap." An inefficient DMU-o is characterized by a positive virtual gap, indicating a discrepancy between its current performance and the Efficiency Equator. The DMU must undergo a comprehensive four-phase improvement process.

\indent The commencement of this process involves the identification of the first SIC scalar, $\kappa_o^1$. This scalar is conceptualized as the sum of the estimated intensities associated with DMU-o's efficient counterparts, $\kappa_o^1 = \sum_{\forall j\in \mathcal{E}_o^{PT}}\pi_{jo}^\star$. Essentially, it serves as a quantitative measure of the extent to which DMU-o needs to adjust its operations to align with the practices of its reference peers.

\indent Conversely, a DMU identified as efficient exhibits a zero virtual gap, signifying that its performance already aligns with the best practices observed across the dataset at the optimal level. Such DMUs are deemed ineligible for further analysis within this framework. For these units, the second scenario presented in Section \ref{sec:3} should be employed for continuous performance management and improvement, ensuring they maintain their efficiency status over time.

\subsection {Phase 2: TS1 Model Uses the First Scalar.}\label{sec:4.2}
The comparison between the optimal solutions of the PT and TS1 models in Steps I and II  reveals some critical insights. Both models share identical TAP solutions in Step I due to the equality $\tau_o^{PT\#} = \tau_o^{TS1\#}$ and $\sum_{\forall j\in \mathcal{E}_o^{PT}} \pi_{jo}^\#  $= $\sum_{\forall j\in \mathcal{E}_o^{TS1} } \pi_{jo}^\#$ = $\kappa_o^1$. However, their TVG programs yield two sets of optimal solutions due to the additional vScalar. This process leads to a specific relationship between their intensities, resulting in disparate benchmarks for performance indices between the two models.

\begin{theorem}\label{thm:10}{ Step I of PT and TS1 models are linked.}
\end{theorem}
\noindent\textbf{Proof.} (\ref{Eq:90}) shows Step I of the PT and TS1 models have the relationships.
\begin{equation} \label{Eq:90}
\begin{aligned}
	&\Delta_o^{TS1\#} = \Delta_o^{PT\#}= \delta_o^{TS1\#} = \delta_o^{PT\#} ; (Q_o^{PT\#}, P_o ^{PT\#})=(Q_o^{TS1\#}, P_o^{TS1\#});\\
	& (v_o^{PT\#}, u_o^{PT\#} )\neq(v_o^{TS1\#}, u_o^{TS1\#}); \sum_{\forall j\in \mathcal{E}_o^{PT}} \pi_{jo}^\# = \sum_{\forall j\in \mathcal{E}_o^{TS1}} \pi_{jo}^\# = \kappa_o^1.   \end{aligned} \end{equation}	
    
 Note the best peers’ intensities, $\pi_{jo}^{PT\#}$, and $\pi_{jo}^{TS1\#}$, may be different; therefore, the benchmarks of performance indices of the two models are distant. They have identical solutions on the dual variables. But have different primal solutions. 
 \hfill $\square$

\indent Step II further differentiates the two models. Although their TAP programs provide the same total values, their goal prices and other variables differ. This disparity leads to different efficiency measurements: $E_o^{PT\star}$  represents pure technical efficiency, and $E_o^{TS1\star}$ comprises the effects on the gray virtual input and output and scalar prices. Notably, the SIC scalar proves to be a crucial factor in evaluating DMU-o. In Step II of PT and TS1 models,  the TAP programs equal to
($\sum_{\forall i\in I}$ $Q_{io}^\star$ 
+ $\sum_{\forall r\in R}$ $P_{ro}^\star$ ) multiplies the distinct virtual goal prices, $\tau_o^{PT\star}$  and $\tau_o^{TS1\star}$. The PT and TS1 models have the relationships shown in  (\ref{Eq:91}).
\begin{equation}\label{Eq:91}
\begin{aligned}
	&\tau_o^{TS1\star} \neq \tau_o ^{PT\star}, \Delta_o^{TS1\star} \neq \Delta_o^{PT\star};\quad(v_o^{PT\star}, u_o^{PT\star}) \neq (v_o^{TS1\star}, u_o ^{TS1});\\ 
	&(Q_o^{PT\star}, P_o ^{PT\star}) = (Q_o^{TS1\star}, P_o ^{TS1\star});\quad 0 < E_o^{PT\star} \neq E_o^{TS1\star} \le1.	
	\end{aligned}
	\end{equation}

\subsection {Phase 3: Identifying the Second Scalar.} \label{sec:4.3} The standard sensitivity analysis of the LP TS1 model involves perturbing the $\kappa_o^1$  scalar, which results in allowable decreases and increases. Specifically, $\kappa_o^2$  equals $\kappa_o^1$  minus the allowable decrease or $\kappa_o^2$  equals $\kappa_o^1$  plus the allowable increase, as expressed in  (\ref{Eq:92}).
\begin{equation} \begin{aligned}
	\kappa_o^2  &= \kappa_o^1  - \textrm{(allowable decreasing of } \kappa_o^1), \\
 \textrm{or  }
	 \kappa_o^2   &= \kappa_o^1   + \textrm{(allowable increasing of }\kappa_o^1 ).
	\end{aligned} 
	\label {Eq:92}
	\end{equation}	
 \vspace{-1 em}
\begin{theorem}\label{thm:11}{The final SIC scalar $\kappa_o^c$ in Phase-4 is linked to $\kappa_o^1$ in Phase-1.}
\end{theorem}
\noindent \textbf{Proof.} These bounds ($\kappa_o^1$  and $\kappa_o^2$) represent the possible range of SIC scalars ($\kappa_o^c $) within the TSc model. If $w_o^{TS1\star}$  and $w_o^{TS2\star}$ are greater than 0, $\kappa_o^1  < \kappa_o^c  < \kappa_o^2 $. Otherwise, $\kappa_o^1  > \kappa_o^c  > \kappa_o^2 $. \hfill $\square$

\indent The best peers in the TS1 and  TS2 models remain consistent, as indicated by $\mathcal{E}_o^{TS1}$=$\mathcal{E}_o^{TS2}$=$\mathcal{E}_o^{TSc}$. However, the estimated intensities of these peers differ, $\pi^{TS1\star}$ $\neq\pi^{TS2\star}$. The solutions of these models exhibit relationships detailed in (\ref{Eq:93}), displaying discrepancies in the estimated variables, such as:
\begin{equation}
\begin{aligned}
 \tau_o ^{TS1\star} \neq \tau_o ^{TS2\star};
\quad(v_o^{TS1\star}, u_o^{TS1\star}) \neq (v_o^{TS2\star}, u_o^{TS2\star});\\  w_o^{TS1\star}\neq w_o^{TS2\star};\quad (Q_o^{TS1\star}, P_o^{TS1\star}) \neq(Q_o^{TS2\star}, P_o^{TS2\star}).
 \label{Eq:93}  
 \end{aligned} 
 \end{equation} 

\subsection{Phase 4: Select the Final Scalar Within the Region.}\label{sec:4.4} 
Phases 1, 2, and 3 integration confirms DMU-o's learning process with the reference peers, whose total intensities are bounded within the $\kappa_o^1$  and $\kappa_o^2$. Phase 4, a critical decision-making step, involves selecting the final scalar $\kappa_o^z$  for the  TSz model. DMU-o undertakes several trials in choosing scalars in Phase 4 to determine a preferred final scalar, $\kappa_o^z$. Each choice scalar corresponds to a unique set of estimated intensities and benchmarks of inputs and outputs,  $\widehat{x}_{io}^{TSz}$  and  $\widehat{y}_{ro}^{TSz}$ derived from the efficient peers.

\indent  By emulating the reference peers based on the estimated intensities, DMU-o gains insights into feasible and desirable productivity benchmarks. This interactive approach ensures a more thorough evaluation of DMU-o. Exploring different SIC scalars and model variations aids in comprehending the diverse impacts these factors have on assessments. Ultimately, this process helps make better-informed management decisions regarding productivity and efficiency.

\begin{theorem}\label{thm:12}{ DMU-o chooses the final SIC scalar $\kappa_o^z$ to fulfill Condition C4, the estimated AARs are applicable.}
\end{theorem}

\noindent\textbf{Proof.} DMU-o may redefine the datasets by altering the performance indices and DMUs and repeating the four-phase procedure. \hfill $\square$

\section{Numerical examples and the 2D geometric intuitions}\label{sec:5} DMU-o (DMU-K) evaluates its performance relative to compatible peers such as DMUs K, A, B, D, G, and H. It selects inputs x1 and x2 and outputs y1 and y2, measured in various units. Some market prices per unit of the inputs and outputs are not attainable. This assessment disregards the interactions between the four indices. DMU-o supplies the dataset detailed in Tables \ref{table:1}, which other DMUs may accept and use for their evaluations. 
     \begin{table}[H]
     \centering
     \renewcommand{\arraystretch}{1.0}
     \caption{The example data.} \label{table:1}
     \begin{tabular}{ccccccc} 
     \hline
DMU-j & K & A & B & D & G & H\\ 
\hline
 $x_{1j}$(ton)   &1.6   &2.3   &1   &1.9   &1.8   &2.5   \\
$x_{2j}$(hr)   &145   &120   &29   &281   &250   &100   \\
\hline
$y_{1j}$$(m^3)   $&1036   &1327   &567   &2446   &1794   &1000   \\
$y_{2j}$($\%$)   &49   &97   &89   &97   &57   &70   \\
\hline
 \end{tabular}
\end{table} 

\indent Table \ref{table:2} summarizes the inefficiency measurements of DMU-o within the first scenario by the PT, TS1, TS2, and TS3 models. Model PT evaluates DMU-B and DMU-D and confirms they are efficient units. Section \ref{sec:5.2} offers the second scenario to measure the super-efficiencies of efficient DMUs.
\subsection{Measures the Inefficiencies}\label{sec:5.1} 
\subsubsection{Relationships Between PT and TS1 Models in Step I}\label{sec:5.1.1}
(\ref{Eq:18}) and (\ref{Eq:35}) are used to calculate the total adjustment prices ($\delta_o^{PT\#}$ = $\delta_o^{TS1\#}$) and total virtual gaps ($\Delta_o^{PT\#}$= $\Delta_o^{TS1\#}$) using specific data provided in row R2 of Table \ref{table:2}. Each amount is \$2.3010. Additionally, the TS1-I model incorporates the SIC impact through the $\textit{vScalar}$, $\omega_o^{1\#}$(=\$2.479), in calculating the total virtual gap. This result demonstrates how the virtual gap, considering the SIC and vScalar in the TS1-I model, maintains the same value as in the PT-I model, indicating a consistent evaluation despite including additional factors.

  \begin{figure}[H]
\centering
\includegraphics[width=0.65\linewidth, height=9cm]{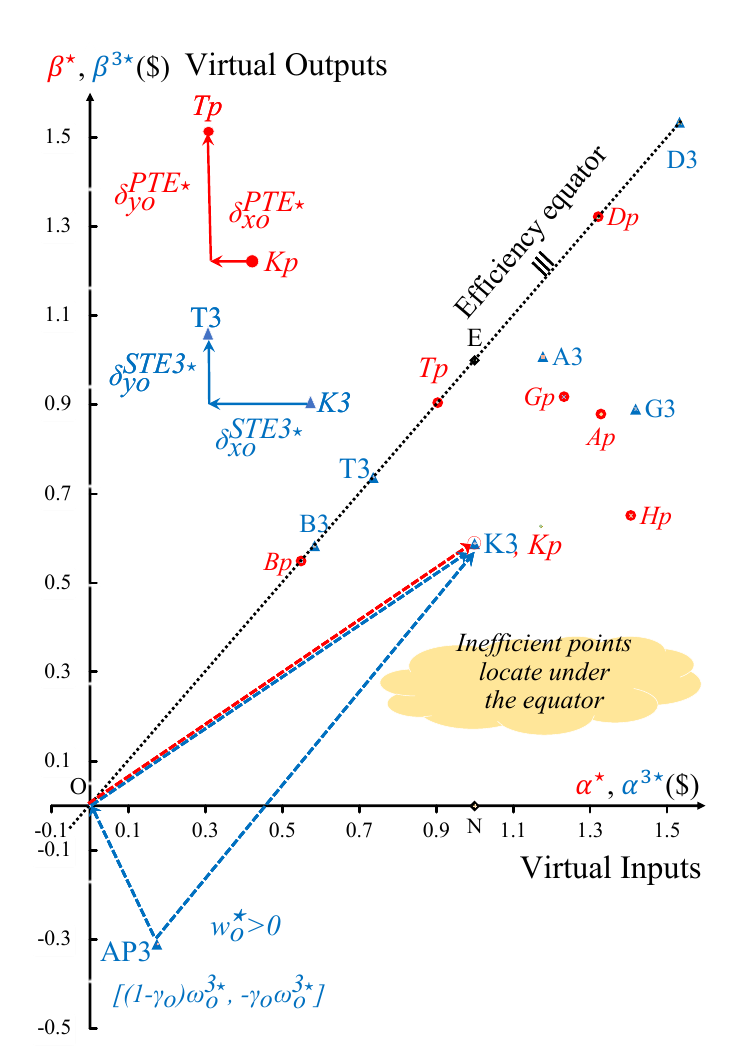} 
\caption{Evaluate DMU-K by PT and TS3 models.}
\label{fig:2}
\end{figure}
 \subsubsection{Compare the PT and TS3 Models}\label{sec:5.1.2}
By try and error, we identified $\kappa_o^3$ equals 0.718 to have $E_o^{PT\star}$=$E_o^{TS3\star}$(=0.589). The columns of PT-II and TS3-II in Table \ref{table:2} list the estimated solutions, in which rows R8 and R9 are the two virtual technology sets, ($\alpha_j^\star,\beta_j^\star$) and ($\alpha_j^{3\star},\beta_j^{3\star}$), where DMU-j= (A, B, D, G, H, K), as per  (\ref{Eq:12}),  (\ref{Eq:28}), and (\ref{Eq:29}).

\indent Via  (\ref{Eq:43}) and (\ref{Eq:44}), the locations of Tp, ($\widehat\alpha_o^\star,\widehat\beta_o^\star$) and point T3, ($\widehat\alpha_o^{3\star},\widehat\beta_o^{3\star}$) on the Efficiency Equator are obtained. The estimated rectilinear distances from points Kp to Tp and from point K3 to T3 are $\delta_o^{PT\star}$, and $\delta_o^{TS3\star}$, see the expressions in the upper-left corner of Figure \ref{fig:2}. Section \ref{sec:2.5.4} illustrated the SIC affected the estimations as shown in (\ref{Eq:52}). 
 \indent In Figure \ref{fig:2}, each point signifies a DMU-j within the virtual technology sets of the PT-II ($\Phi_o^{PT\star}$) and TS3-II models ($\Phi_o^{TS3\star}$), respectively. 
 The points (Ap, Bp, Dp, Gp, Hp, Kp) and (A3, B3, D3, G3, H3, K3) located at ($\alpha_j^\star,\beta_j^\star$) and ($\alpha_j^{3\star},\beta_ j^{3\star}$). The blue dashed lines are graphically representing the vectors $\overrightarrow{\rm O, K3}$, $\overrightarrow{\rm O, AP3}$, and $\overrightarrow{\rm AP3, K3}$. The red dashed line indicates the vector $\overrightarrow{\rm O, Kp}$ of pure technical efficiency,  $E_o^{PT\star}$= $\beta_o^{3\star}$/$\alpha_o^{3\star}$. 
\subsubsection {Compare the Solutions between TS1 and TS2 models} \label{sec:5.1.3} 
\begin{figure}[H]
\centering
\includegraphics[width=0.65\linewidth, height=9cm]{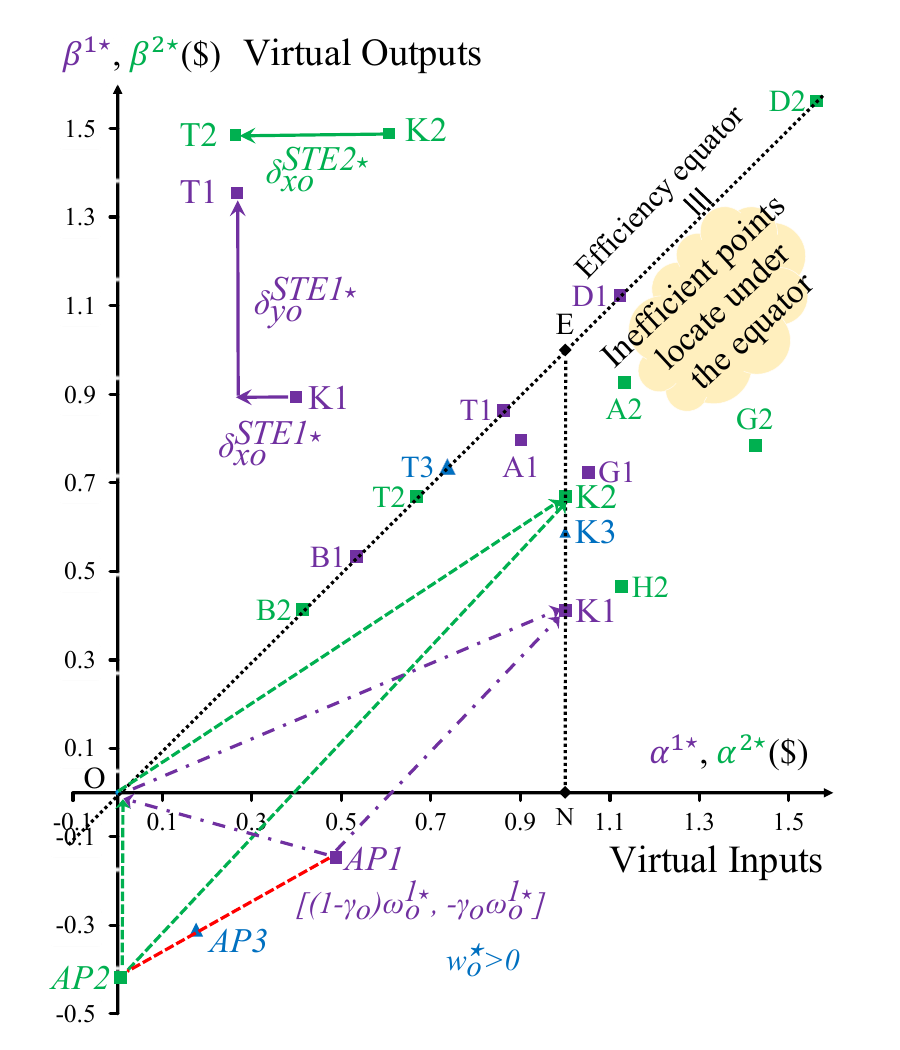}
\caption{Evaluate DMU-K by TS1 and TS2 models. \hfill}
\label{fig:3}
\end{figure}

Figure \ref{fig:3}, similar to Figure \ref{fig:2}, depicts the locations of DMUs in the TS1-II and TS2-II models, respectively, (A1, B1, D1, G1, H1, K1)  and (A2, B2, D2, G2, H2, K2). Their coordinates, ($\alpha_j^{1\star},\beta_j^{1\star}$) and ($\alpha_j^{2\star},\beta_j^{2\star}$), are listed in rows R8 and R9.  DMUs B and D emerge as efficient peers with their intensities $\pi_B^\star$  and $\pi_D^\star$.  R4 detailed the benchmarks of inputs and outputs. R5 lists the virtual prices. R6 contains the benchmark virtual scales. R7 lists $\Xi_o^{TSc\star}$ and $E_o ^{TSc\star}$. Since $w_o^\star>0$, decreasing $\kappa_o^1$  (=1.5153)  to  $\kappa_o^2$  (=0.515), DMU-o exhibits decreasing RTvS: $\Xi_o^{TS1\star} (=2.435) > \Xi_o^{TS2\star} (=1.497)$.
\indent At the same time, the relative efficiency is increasing, $E_o^{TS1\star}$ (=0.411) less than $E_o^{TS2\star}$ (=0.688). However, DMU-o may use  (\ref{Eq:48}) to analyze the interlinkage relationships between input and output indices for managing the performance indices.

  The TSc program introduces an \textit{anchor point}  (APc) located at the point $[(1 -\gamma_o)\kappa_o^c  w_o^\star$, $- \gamma_o \kappa_o^c w_o^\star]$, denotes as [AP(x), AP(y)]. Anchor points are located in the second and fourth quadrants when $w_o^\star> 0$ and $w_o^\star< 0$.
 The rectilinear distance between point APc and the origin O is the vScalar, $\kappa_o^c w_o^\star$. In Figure \ref{fig:3}, the anchor points of the TS1 and TS2 models are located at  AP1 and  AP2. Because $w_o^\star>0$, the two anchor points are located in the fourth quadrant.

 \quad The TS3 model interpolates TS1 and TS2 models, providing a framework for determining the optimal scalars $\kappa_o^3$ within the interval of SIC scalars. Using the scalar $\kappa_o^3$ between $\kappa_o^1$  and $\kappa_o^2$, the TS3 model has the anchor points, and AP3 is located between AP1 and AP2. Inefficient DMUs are located under the Efficiency Equator. 
 \quad  DMU-j lies in the fourth quadrant if, according to  (\ref{Eq:2}), one could have $\Delta_j^{TSc\star}> \$0$. ($\alpha_j^{c\star}>\$0, \beta_j^{c\star}<\$0$), and their ratio is less than zero. Conversely, DMU-j is in the third quadrant if ($\alpha_j^{c\star}<\$0, \beta_j^{c\star}<\$0$), and their ratio is greater than zero. 

\quad Section \ref{sec:2.5.4} illustrated the SIC affected the estimations as shown in (\ref{Eq:53}). Therefore, eliminating an inefficient DMU-j from (\textit{X, Y}) does not impact the evaluation. Using geometric vectors, we can visually interpret the effects under the SIC scalars $\kappa_o^1$  and $\kappa_o^2$. In Figure \ref{fig:3}, the purple color dashed lines express the three vectors $\overrightarrow{\rm O,K1},  \overrightarrow{\rm AP1,O}$, and $\overrightarrow{\rm AP1,K1}$ of the TS1 model can be found.

 \indent  Similarly, the green color dashed lines express the three vectors $\overrightarrow{\rm O, K2},  \overrightarrow{\rm AP2, O}$, and $\overrightarrow{\rm AP2, K2}$ of the TS2 model. The three vectors $\overrightarrow{\rm O,K3}$, $\overrightarrow{\rm AP3,O}$ and $\overrightarrow{\rm AP3, K3}$ of TS3 model are not depicted for simplicity.

\begin{table}[H]
\centering
\setlength{\tabcolsep}{3pt}
\renewcommand{\arraystretch}{0.91}
\caption{Evaluate Inefficiencies for DMU-K, DMU-B, and DMU-D.} \label{table:2}
\footnotesize
\begin{tabularx}{\linewidth}{lccrrrrrrrrr} 
  & 'DMU-o=' &  & K & K & K & K & K & K & K & B & D\\ 
 \hline
  Row& Solution & Unit &PT I&PT II &TS1 I&TS1 II&TS2 I&TS2 II &TS3 II&PT II&PT II\\
  \hline

R1 &$\kappa_o^c$   & - &1.5153 &1.5153 &1.5153 &1.5153 &0.5150 &0.5150& 0.718&1 &1\\
   &$\tau_o$ &  $\$$ &1 &0.179 &1.000 &0.256 &1.000 &0.500 &0.413&0.500&	0.266\\
\hline
R2 &$\Delta_o,\delta_o$ & \$ &2.3010 &0.4113 &2.3010 &0.5893 &0.6643 &0.3321&0.411 &0 &0 \\
   &$v_{1o}$  &  \$/ton &2.8713 &0.5133 &0.6250 &0.1601 &0.6250 &0.3125 &0.258&0.500	&0.387 \\
   &$v_{2o} $  &  \$/hr &0.0069 &0.0012 &0.0069 &0.0018 &0.0069 &0.0034&0.003 & 0.017	&0.001\\
   &$u_{1o}$ &  \$/$m^3$ &0.0022 &0.0004 &0.0011 &0.0003 &0.0011 &0.0006&0.0005&0.001&	0.0003 \\
   &$u_{2o}$  &  \$ / $\%$ &0.0204 &0.0036 &0.0204 &0.0052 &0.0204 &0.0102 &0.008&0.006&	0.003 \\
   &$w_o$ &  \$  &0 &0 &1.6362 &0.4190 &1.6362 &0.8181&0.675 &0&	0 \\
\hline
R3 &$Q_{1o} $  & - &0 &0 &0 &0 &0.4554 &0.4554&0.363&0 &0 \\
   &$Q_{2o}$   &-  &0.5334 &0.5334 &0.5334 &0.5334 &0.2089 &0.2089 &0.275&0 &0 \\
   &$P_{1o}$   &-  &0 &0 &0 &0 &0 &0&0&0 &0 \\
   &$P_{2o} $  &-&1.7677 &1.7677 &1.7677 &1.7677 &0 &0 &0.359&0 &0 \\
   & $\pi_B$ &- &1.421 &1.421 &1.421 &1.421 &0.119 &0.119&0.383&1&0 \\
   & $\pi_D$  &- &0.094 &0.094 &0.094 &0.094 &0.396 &0.396 &0.335&0 &1 \\
   \hline
R4 &$\widehat{x}_{1o}$ &ton &1.6 &1.6 &1.6 &1.6 &0.8713 &0.8713&1.019 &1.0&	1.9\\
   &$\widehat{x}_{2o}$  &hr &67.66 &67.66 &67.66 &67.66 &114.72 &114.72 &105.165&29	&281\\
   & $\widehat{y}_{1o}$  & $m^3$ &1036 &1036 &1036 &1036 &1036 &1036 &1036&567&	2446 \\
   & $\widehat{y}_{2o}$ &\%  &135.6 &135.6 &135.6 &135.6 &49.0 &49.0&66.58&89&	97\\
   & $\gamma_o$ & -  &0.232 &0.232 &0.232 &0.232 &1.000 &1.000&0.640&0&0 \\
   &  $\omega_o^a$ &  \$ &0 &0 &2.479 &0.635 &0.843 &0.421&0.485 &0 &0 \\
   \hline
R5 &$x_{1o}v_{1o}$ &  \$ &4.594 &0.821 &1.000 &0.256 &1.578 &0.789&0.689&0.500&	0.734 \\
   &$x_{2o}v_{2o}$ &  \$ &1.000 &0.179 &3.479 &0.891 &1.265 &0.632 &0.622&0.500&	0.266\\
   &$y_{1o}u_{1o}$ &  \$ &2.293 &0.410 &1.178 &0.302 &1.178 &0.589 &0.486&0.500	&0.734 \\
   &$y_{1o}u_{1o}$ &  \$ &1.000 &0.179 &3.479 &0.891 &1.843 &0.921 &0.898&0.500	&0.266 \\
\hline
R6 &$\alpha_o,\alpha_{o}^c$&  \$ &5.594 &1.000 &3.905 &1.000 &2.000 &1.000&1.000 &1 &1 \\
&$\beta_o,\beta_{o}^c$ &  \$ &3.293 &0.589 &1.604 &0.411 &1.336 &0.668&	0.589 &1 &1 \\
&$\widehat\alpha_o,\widehat\alpha_{o}^c$ &  \$ &5.061 &0.905 &3.371 &0.863 &1.336 &0.668&0.737 &1&1 \\
   &$\widehat\beta_o,\widehat\beta_{o}^c$ &  \$&5.061 &0.905 &3.371 &0.863 &1.336 &0.668 &0.737&1 &1\\
\hline
R7 & $\Xi_o$ &- &1.699 &1.699 &2.435 &2.435 &1.497 &1.497 &1.699&1&1 \\
   & $E_o$ & - &0.589 &0.589 &0.411 &0.411 &0.668 &0.668&0.589 &1 &1 \\
   \hline 
   R8&$\alpha_K,\alpha_{K}^c$&  \$&5.594&	1&	3.905	&1.000&	2&	1&	1.000&	3.300	&0.755\\
  &$\alpha_A,\alpha_{A}^c$&  \$ &7.432 &1.000 &3.522 &0.902 &2.265 &1.133&0.902&3.219	&1.002 \\ 
   &$\alpha_B,\alpha_{B}^c$&  \$ &3.071 &1.328 &2.082 &0.533 &0.825 &0.413 &0.533&1&	0.414 \\   
      &$\alpha_D,\alpha_{D}^c$&  \$ &7.393 &0.549 &4.382&1.122 &3.125 &1.563 &1.122&5.795&	1 \\
   & $\alpha_G,\alpha_{G}^c$ &  \$ &6.892 &1.322 &4.106 &1.052 &2.849 &1.425 &1.052&5.210&	0.932 \\
   & $\alpha_H,\alpha_{H}^c$&  \$ &7.868 &1.232 &3.509 &0.899 &2.252 &1.126 &0.899&2.974&	1.061 \\  
   \hline 
   R9&$\beta_K,\beta_{K}^c$&  \$&3.293&	0.589&	1.604&	0.411&	1.336&	0.668&	0.589&	1.189&	0.445\\
&$\beta_A,\beta_{A}^c$ &  \$ &4.917 &0.879 &3.110 &0.796 &1.853 &0.926 &0.796&1.715&	0.664 \\
   &$\beta_B,\beta_{B}^c$ &  \$ &3.071 &0.549 &2.082 &0.533 &0.825 &0.413 &0.533&1&	0.414 \\
   &$\beta_D,\beta_{D}^c$ &  \$ &7.393 &1.322 &4.382 &1.122 &3.125 &1.563 &1.122&2.702&	1 \\
   &$\beta_G,\beta_{G}^c$&  \$ &5.134 &0.918 &2.824 &0.723 &1.568 &0.784&0.723 &1.902	&0.695 \\
   &$\beta_H,\beta_{H}^c$ &  \$&3.642 &0.651 &2.187 &0.560 &0.930 &0.465 &0.560&1.275	&0.492\\
   \hline
  R10&AP(x)&\$& -&	-	&1.905	&0.488&	0&	0&	0.175&	-	&-\\
&AP(y)&\$&-&	-&	-0.575	&-0.147	&-0.843&	-0.421	&-0.310	&-	&-\\ \hline
 \end{tabularx}
 \end{table}

\indent Point K3 is located between points K1 and K2 on the normalization vertical line (N, E).  Using the final SIC scalar $\kappa_o^z$ in the TSz model, points APz, Kz, and Tz are located between points (AP1 and AP2), (K1 and K2), and (T1 and T2), respectively. The three vectors, $\overrightarrow{\rm O,Kz}$, $\overrightarrow{\rm APz,O}$, $\overrightarrow{\rm APz,Kz}$, shall be visualized in Figure \ref{fig:3}, similar to (\ref{Eq:53}). 
The upper left corner of the figure depicts the rectilinear distances of DMU-o to project on the Efficiency Equator.
\subsection{Measures the Super-efficiencies}\label{sec:5.2} 
\begin{figure}[H]
\centering
\includegraphics[width=0.63\linewidth, height=8.5cm]{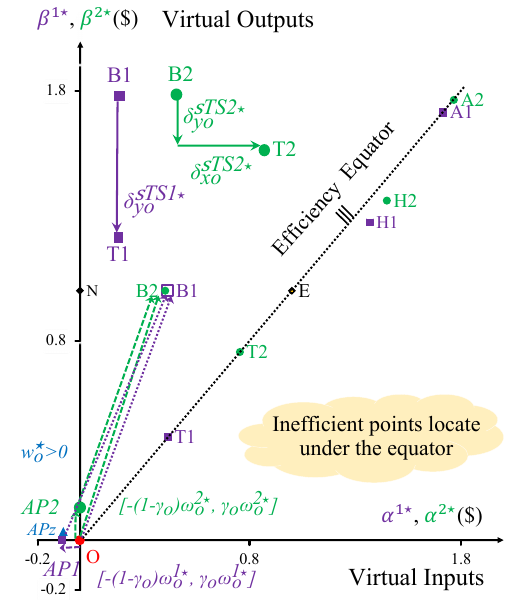} 
\caption{Evaluate DMU-B by sTS1 and sTS2 Models.}
\label{fig:4}
\end{figure}

The first scenario confirmed DMU-B and DMU-D are efficient. The second scenario employs the sPT and sTSc models to measure their super-efficiencies in the four-\noindent phase procedure. Table \ref{table:3} summarized the solutions of the sPT, sTS1, sTS2, and sTSz models. To illustrate the fourth phase, we assume the final SIC scalars are the mid-point of the interval of $\kappa_o^1$ and $\kappa_o^2$. The anchor point of the sTSz model that contains the SIC scalar  $\kappa_o^z$, APz is located in the region of points O, AP1, and AP2 because DMU-o is excluded from (\ref{Eq:59}). In Figure \ref{fig:4}, AP1 and AP2 are in the second quadrant because $\omega_o^c>\$0$. When $\omega_o^c>\$0$, the anchor point will be located in the fourth quadrant.

 \indent R6 in Table \ref{table:3} listed the locations of DMU-o and its projection points on the Efficiency Equator. DMU-o is located in the first quadrant and above the Efficiency Equator. \noindent DMU-o projects on the Efficiency Equator, the T1 and T2 in Figure \ref{fig:4} and \ref{fig:5}.
 The reference DMUs of DMU-B, A1, and A2 are on the Efficiency Equator, while the remaining inefficient DMUs lie under it. In evaluating DMU-D, DMU-B, and DMU-G, they are located on the Efficiency Equator.

\indent In Table \ref{table:3}, DMU-B has higher super-efficiencies (see R7), which reduced the outputs will have an efficiency score equal to 1. Regarding EA, DMU-B has strong outputs and can be selected as the best alternative to the MCDM problem.

\section {Discussions} \label{sec:6}

\indent Applying the four-phase procedure to an MCDM problem, the decision-maker needs to manipulate the SIC scalar $\kappa_o^z$ for each DMU-o. The decision-maker performs the sensitivity analysis to compare the choices. The efficient DMUs may not be selected since the uncertainties could not always be expressed in the data matrix and the EA solutions. Our contributions to the MCDM problems are as follows.
\begin{romanlist}[(ii)]
\item Integration for Enhanced Decision Support: MCDM can incorporate efficiency analysis results as criteria,  providing a quantitative basis for evaluating alternatives. For instance, VGA results can serve as efficiency scores within an MCDM framework to rank and select the most efficient DMUs or alternatives. 
\item Complementary in Performance Improvement: 
While MCDM facilitates the selection among alternatives based on multiple criteria, efficiency analysis identifies gaps and potential areas for deterioration. Integrating these approaches can guide strategic decisions on resource allocation, process optimization, and performance enhancement.
\begin{figure}[H]
\centering
\includegraphics[width=0.63\linewidth, height=8.5cm]{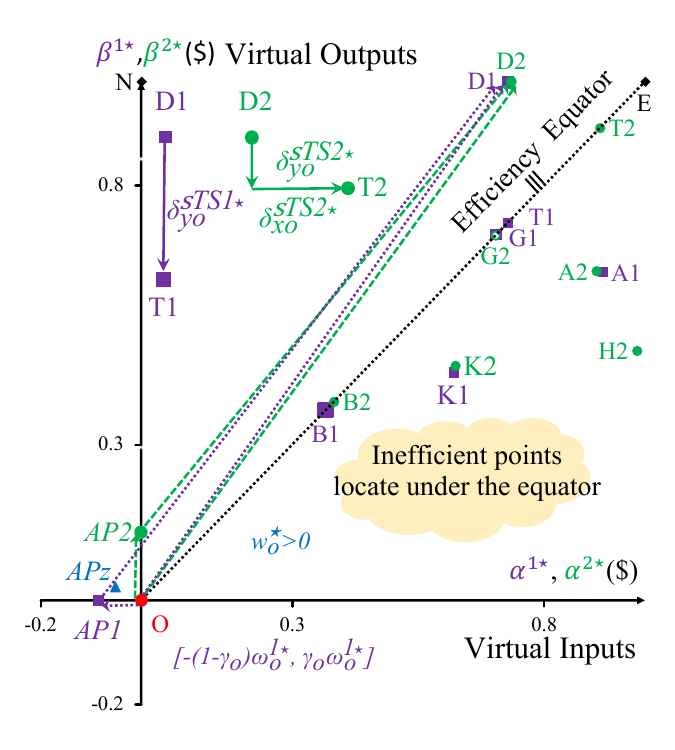}
\caption{Evaluate DMU-D by sTS1 and sTS2 Models.}
\label{fig:5}
\end{figure}
 \item Framework for Complex Decisions: In complex decision 
 scenarios, such as strategic planning, resource management, and system design, combining MCDM with efficiency analysis offers a robust framework for considering a broad spectrum of performance metrics and stakeholder preferences, thereby supporting more informed and comprehensive decisions.
 \item Methodological Cross-Fertilization: Advances in MCDM and efficiency analysis methodologies can enrich each other. For example, MCDM techniques can be adapted to address specific challenges in efficiency analysis, such as handling qualitative conditions or integrating stakeholder preferences into efficiency models.
 \end{romanlist}
 
 \indent Further research could explore developing integrated models that leverage the strengths of both MCDM and EA, tailored to specific sectors such as healthcare, energy, and public services. Methodological innovations in handling uncertainty, dynamic Conditions, and the inclusion of sustainability dimensions are promising  
 \newpage
\begin{table}[t]
\centering
\setlength{\tabcolsep}{3 pt}
\renewcommand{\arraystretch}{0.91}
\caption{Solutions in Measuring the Super-Efficiencies of DMU-B and DMU-D.} 
\footnotesize
\begin{tabularx}{\linewidth}{cccrrrrrrrr} \label{table:3} 
 & 'DMU-o=' & &B &B &B &B &D &D &D &D \\
\hline 
{Row}&{Solution}& Unit &sPT-II&sTS1-II&sTSz-II&sTS2-II&sPT-II&sTS1-II&sTSz-II&sTS2-II\\
  \hline
R1&$\kappa_o^c $ &  -    & 0.2417 & 0.2417 & 0.3345 & 0.4273 & 1.3411 & 1.3411 & 1.4092 & 1.4774 \\
  &$\tau_o $ &   \$    & 0.5   & 0.5   & 0.4823 & 0.4591 & 0.7709 & 0.8037 & 0.7827 & 0.7614 \\
 & $\Delta_o,\delta_o$& \$ & 0.5855 & 0.5855 & 0.5964 & 0.5979 & 0.2611 & 0.2722 & 0.2688 & 0.2652 \\
 \hline
  R2 &$v_{1o} $  & \$/ton & 0     & 0     & 0     & 0     & 0.3889 & 0.4230 & 0.4119 & 0.4007 \\
   &$v_{2o}$ &   \$/hr    & 0.0143 & 0.0172 & 0.0166 & 0.0158 & 0     & 0     & 0     & 0 \\
   &$u_{1o}$ &    \$/ $m^3$    & 0.0009 & 0.0009 & 0.0009 & 0.0008 & 0.0003 & 0.0003 & 0.0003 & 0.0003 \\
  &$u_{2o}$ &   \$/\%    & 0.0056 & 0.0056 & 0.0054 & 0.0052 & 0.0024 & 0.0020 & 0.0020 & 0.0019 \\
&$w_o$    &    \$   & 0.0000 & 0.3538 & 0.3413 & 0.3249 & 0     & 0.0566 & 0.0551 & 0.0536 \\
          \hline
   R3&$Q_{1o}$ &   -    & 0     & 0     & 0     & 0     & 0     & 0.0000 & 0.1156 & 0.2313 \\
   &$Q_{2o}$  &   -    & 0     & 0     & 0.3840 & 0.7681 & 0     & 0     & 0     & 0 \\
   &$P_{1o} $   &   -    & 0.4344 & 0.4344 & 0.2172 & 0.0000 & 0.3387 & 0.3387 & 0.2278 & 0.1170 \\
   &$P_{2o} $  &   -    & 0.7366 & 0.7366 & 0.6355 & 0.5343 & 0     & 0     & 0     & 0 \\
&$\pi_A$  &   -    & 0.2417 & 0.2417 & 0.3345 & 0.4273 & 0     & 0     & 0     & 0 \\
&$\pi_B$ &   -    & 0     & 0     & 0     & 0     & 0.6424 & 0.6424 & 0.5211 & 0.3997 \\
&$\pi_G$   &    -   & 0     & 0     & 0     & 0     & 0.6986 & 0.6986 & 0.8881 & 1.0776 \\
\hline
R4 &$\widehat{x}_{1o}$  &   ton    & 1     & 1     & 1     & 51.273 & 1.9   & 1.9   & 2.119 & 2.339 \\
 &$\widehat{x}_{2o}$    &   hr    & 29    & 29    & 40.137 & 567   & 281   & 281   & 281   & 281 \\
& $\widehat{y}_{1o}$   &   $m^3$    & 320.7 & 320.7 & 443.8 & 41.4  & 1617.6 & 1617.6 & 1888.8 & 2159.9 \\
 & $\widehat{y}_{2o}$ &    \%   & 23.442 & 23.442 & 32.444 & 0.5897 & 97    & 97    & 97    & 97 \\
   & $\gamma_o$ &    -   & 0     & 0     & 0.3105 & 0.1388 & 0     & 0.0000 & 0.3367 & 0.6642 \\
  & $\omega_o^a$ &   \$ & 0     & 0.0855 & 0.1142 & 0     & 0     & 0.0759 & 0.0776 & 0.0792 \\
  \hline
R5 &$x_{1o}v_{1o}$&  \$     & 0     & 0     & 0     & 0.4591 & 0.7389 & 0.8037 & 0.7827 & 0.7614 \\
&$x_{2o}v_{2o}$ & \$   &  0.4145 & 0.5   & 0.4823 & 0.4591 & 0     & 0     & 0     & 0 \\
   &$y_{1o}u_{1o}$  &  \$     & 0.5   & 0.5   & 0.4823 & 0.4591 & 0.7709 & 0.8037 & 0.7827 & 0.7614 \\
 &$y_{2o}u_{2o}$ &   \$    & 0.5   & 0.5   & 0.4823 & 0.4021 & 0.2291 & 0.1963 & 0.1912 & 0.1860 \\
\hline
 R6 &$\alpha_o,\alpha_{o}^c$   &   \$    & 0.4145 & 0.4145 & 0.4036 & 1     & 0.7389 & 0.7278 & 0.7312 & 0.7348 \\
&$\beta_o,\beta_{o}^c$   &   \$    & 1     & 1     & 1     & 0.7547 & 1     & 1     & 1     & 1 \\
 &$\widehat\beta_o,\widehat\beta_{o}^c$  &    \$   & 0.4145 & 0.4145 & 0.5888 & 0.7547 & 0.7389 & 0.7278 & 0.8217 & 0.9109 \\
 &$\widehat\beta_o,\widehat\beta_{o}^c$  &    \$   & 0.4145 & 0.4145 & 0.5888 & 0.4021 & 0.7389 & 0.7278 & 0.8217 & 0.9109 \\
 \hline
 R7   & $\Xi_o$  &   -    & 0.4145 & 0.4145 & 0.4036 & 2.4868 & 0.7389 & 0.7278 & 0.7312 & 0.7348 \\
 & $E_o$ &   -    & 2.4126 & 2.4126 & 2.4779 & 2.1621 & 1.3533 & 1.3740 & 1.3676 & 1.3609 \\
 \hline
R8& $\alpha_K,\alpha_{K}^c$ &  \$     & 2.0725 & 2.1462 & 2.1761 & 1.7663 & 0.6223 & 0.6202 & 0.6226 & 0.6232 \\
&$\alpha_A,\alpha_{A}^c$ &   \$    & 1.7151 & 1.7151 & 1.7603 & 0.4021 & 0.8945 & 0.9163 & 0.9109 & 0.9037 \\
&$\alpha_B,\alpha_{B}^c$ &   \$    & 0.4145 & 0.4145 & 0.4036 & 4.3150 & 0.3889 & 0.3664 & 0.3754 & 0.3827 \\
&$\alpha_D,\alpha_{D}^c$&    \$   & 4.0163 & 4.4910 & 4.4378 & 3.8242 & 0.7389 & 0.7278 & 0.7312 & 0.7348 \\
&$\alpha_G,\alpha_{G}^c$&     \$  & 3.5732 & 3.9565 & 3.9222 & 1.4497 & 0.7000 & 0.7048 & 0.7049 & 0.7033 \\
&$\alpha_H,\alpha_{H}^c$&   \$    & 1.4293 & 1.3703 & 1.4277 & 1.2831 & 0.9723 & 1.0009 & 0.9933 & 0.9839 \\
  \hline
R9 & $\beta_K,\beta_{K}^c$   &   \$    & 1.1889 & 1.1889 & 1.2527 & 1.7663 & 0.4422 & 0.4396 & 0.4466 & 0.4520 \\
 & $\beta_A,\beta_{A}^c$ &    \$   & 1.7151 & 1.7151 & 1.7603 & 1     & 0.6473 & 0.6323 & 0.6343 & 0.6347 \\
 & $\beta_B,\beta_{B}^c$  &    \$   & 1     & 1     & 1     & 2.6723 & 0.3889 & 0.3664 & 0.3754 & 0.3827 \\
& $\beta_D,\beta_{D}^c$  &   \$    & 2.702 & 2.702 & 2.712 & 1.938 & 1     & 1     & 1     & 1 \\
& $\beta_G,\beta_{G}^c$   &   \$    & 1.902 & 1.902 & 1.941 & 1.362 & 0.7000 & 0.7048 & 0.7049 & 0.7033 \\
 & $\beta_H,\beta_{H}^c$   &   \$    & 1.2751 & 1.2751 & 1.3359 & 0     & 0.4805 & 0.4702 & 0.4765 & 0.4811 \\
          \hline
R10&AP(x)&    \$   & 0     & -0.0855 & -0.0787 & 0.0000 & 0     & -0.076 & -0.051 & -0.027 \\
&AP(y)    &  \$     & 0     & 0     & 0.0354 & 0.1313 & 0     & 0     & 0.026 & 0.053 \\ \hline
\end{tabularx}
\end{table}
\clearpage

 \noindent areas that can enhance the applicability and impact of these decision-support tools. In conclusion, the interplay between MCDM methods and EA represents fertile ground for advancing sophisticated  and practical decision-making frameworks, offering significant potential for improving system performance  
 and achieving strategic goals across diverse contexts.

\indent The most effective way to help the decision-making unit become the best performer is by conducting the structured four-phase procedure depicted in Figure \ref{fig:1} that accurately estimates the necessary and achievable input and output adjustments. This procedure is iterative and may require refinement as new data become available or the operational context of the DMUs changes.

 Theorems \ref{thm:1} to \ref{thm:12} have proved that DMU-o is comprehensively assessed. DMU-o will mimic the best peers. Collaboration with domain experts and stakeholders is crucial throughout the process to ensure the relevance and applicability of the findings and recommendations.

\indent The Simplex method [Dantzig and Thapa]\cite{9} is the backbone of the VGA models, enabling post-optimality analysis and flexibility in evaluating DMU-o. This analysis allows for adjustments based on different coefficients $\textit(X, Y)$, which is crucial in real-world evaluations. 

\indent The adjustment ratios of inputs ($Q_o^\star$) typically range between 0 and 1, but outputs ($P_o^\star$) may exceed 1. Imposing limits like $P_r \le 1$  or setting upper and lower bounds, $\underline{Q}_o\le Q_o \le \bar{Q}_o$ and $\underline{P}_o\le P_o \le \bar{P}_o$, to meet its particular requirements. 

\indent Panwar et al.\cite{10} reviewed the practical issues that have been addressed in the DEA literature, such as the \textit{Malmquist Index}, \textit{network DEA}, \textit{free-disposal hull}, \textit{imprecise data}, \textit{data mixed with 0s and negatives} [Sueyoshi and Goto]\cite{11} input/output sharing, and input/output selection [Ali and Seiford]\cite{12}. The current best practice of the four-phase VGA procedure could be used in specific studies to evaluate performance under the worst practice.

\section* {Acknowledgment} 
We thank the referees for their valuable comments. This research received no specific grant from public, commercial, or not-for-profit funding agencies. We have no competing interests to declare. We confirm that the data supporting the findings of this study are available within the article.

\appendix
\section{CRS and VRS DEA Models} \label{appendix: AAAAA}
\subsection{DEA Theories.} Charnes et al.\cite{13} introduced the concept of Data Envelopment Analysis (DEA), which evaluates a set of entities known as decision-making units (DMUs). Each DMU operates a system that consumes multiple inputs to produce outputs. The observed data from these DMUs form a \textit{production possibility set (PPS)} with an efficiency frontier, where each DMU-o $(x_o, y_o)$ belongs to the PPS. The PPS possesses four properties, denoted as (A1)–(A4)[Chapter 3, Cooper et al.]\cite{4}. Banker et al.\cite{14} and Banker et al.\cite{15} introduced the \textit{convexity condition} to the PPS. Additionally, Chapters 1 and 2 in Sickles and Zelenyuk\cite{5} provided twenty graphical illustrations of the PPS in two dimensions to underpin the theoretical foundations of CRS and VRS DEA models. 
 
\indent Eqs. (3.1) and (4.1) in Cooper et al.\cite{4}, as shown 
 in Eqs. (\ref{A.1}) and (\ref{A.2}) represent the PPS without and with the convexity condition, respectively, under the constant returns-to-scale (CRS) and variable returns-to-scale (VRS) assumptions.
\begin{equation} \label{A.1}
    P^{CRS}=\{(x_o,y_o)|x_o\leq X\lambda, y_o \geq Y\lambda, \lambda\geq 0.\}
\end{equation}
\vspace{-2em}
\begin{equation} \label{A.2}
    P^{VRS}=\{(x_o,y_o)|x_o\leq X\lambda, y_o \geq Y\lambda, (e\lambda=1,) \lambda\geq 0.\}
\end{equation}

\indent Cooper et al.\cite{4} offered the envelopment and multiplier programs of the Additive Model, Eqs. (4.34)–(4.38) and Eqs. (4.39)–(4.43), to estimate the vector of variable intensities $\lambda$. Halická and Trnovská\cite{16} extensively reviewed DEA models and proposed a unified model to encompass various DEA approaches that address the Additive Model.
 \subsection {Use the Additive DEA Model to Solve the CRS and VRS PPS.} In the following Eqs. (\ref{A.3})–(\ref{A.7}) and Eqs. (\ref{A.8})–(\ref{A.11}) depicts the envelopment and multiplier program of the unified DEA model. Including and excluding the parenthesized elements, the program represents the VRS and CRS measurement. In which, Eqs. (\ref{A.4})–(\ref{A.7}) constitute Eqs. (\ref{A.1}) and (\ref{A.2}).
 
The \textit{Additive Model} measures the inefficiency of the object DMU, DMU-o, to obtain efficient peers in its regards. Every DMU in the set of DMUs takes turns to be DMU-o to constitute the entire efficiency frontier of the PPS. 
\subsubsection{The Envelopment (primal) Program.}
\begin{equation}
	F_o^\star= \max_{\lambda_o\textbf, \textbf{s}_o^x, \textbf{s}_o^y} \sum_{\forall i\in I} s_{io}^x b_i^x  + \sum_{\forall r\in R} s_{ro}^y b_r^y, \forall o \in J;	
 \label{A.3} \end{equation}
 \vspace{-1.5 em}
 \begin{equation}
	s.t. \sum_{\forall j\in J } x_{ij} \lambda_{oj} + s_{io}^x   = x_{io}, \forall i \in I; \label{A.4}	
 \end{equation}
 \vspace{- 1.5 em}
\begin{equation}
	- \sum_{\forall j\in J} y_{rj} \lambda_{oj} + s_{ro}^y  = - y_{ro}  , \forall r \in R; \label{A.5}	 
  \end{equation}
  \vspace{-1.5 em}
 \begin{equation}
	(\sum_{\forall j\in J }  \lambda_{oj} = 1;) \label{A.6}	
   \end{equation}
   \vspace{-1.5 em}
\begin{equation}
	\lambda_{oj},\forall j \in J; s_{io}^x, \forall i \in I; s_{ro}^y, \forall r \in R  \ge 0. 	\label{A.7}	
  \end{equation}

\indent Each DEA model uses the variable slacks of input-i and output-r, $s_{io}^x$ and $s_{ro}^y$, to ensure that DMU-o projects on the efficiency frontiers of the two PPSs. The coefficients $b_i^x$ and $b_r^y$ are artificial goal weights corresponding to input-i and output-r. Halická and Trnovská\cite{16} listed objective functions and artificial goal weights used in several DEA models in  Tables 2 and  3.
\subsubsection{The Multiplier (dual) Program.}
The envelopment program can be converted to the multiplier program, Eqs. (\ref{A.8})–(\ref{A.11}).  Usually, the multiplier program is not presented when introducing a DEA model, and the duality characteristics of the model are not verified. As a result, the solutions of the multiplier model are not examined.
\begin{equation} \begin{aligned}
	f_o^\star = \min_{\textbf{v}_o, \textbf{u}_o, \sigma_o}
	 \sum_{\forall i\in I} v_{io} x_{io}  - \sum_{\forall r\in R} u_{ro} y_{ro} (+ 1 \sigma_o), \forall o \in J;
\end{aligned} \label{A.8} 
\end{equation}
\vspace{-1em}
  \begin{equation}
s.t. \sum_{\forall i\in I} v_{io} x_{ij} - \sum_{\forall r\in R} u_{ro} y_{rj} ( + 1 \sigma_o ) \ge 0, \forall j \in J; \label{A.9}	
\end{equation}
\vspace{- 1 em}
 \begin{equation}
v_{io} \ge b_i^x,\forall i \in I; 
u_{ro} \ge b_r^y, \forall r \in R;	\label{A.10}	
\end{equation}
\vspace{- 1 em}
\begin{equation}
v_{io}, \forall i \in I; u_{ro}, \forall r \in R(, \sigma_o) free. \label{A.11}	
 \end{equation}

\indent The \textit{inefficiency condition} Eq. (\ref{A.9}) of each DMU-j bounds its \textit{excess} at a minimum of zero and corresponds to a \textit{variable intensity} $\lambda_{oj}$ in the primal program. The model identifies a set of efficient \textit{reference peers} with zero inefficiency score and positive intensity, while other DMUs are deemed inefficient with zero intensity. 
 The \textit{weighting condition} Eq. (\ref{A.10}) of an input (or output) corresponds to a \textit{variable slack},  $s_{io}^x$ (or $s_{ro}^y$), in the primal program, and restricts its \textit{variable weight}, $v_{io}$ (or $u_{ro}$), to a lower bound defined by a specified \textit{artificial goal weight}, $b_i^x$ (or $b_r^y$).
 \subsection{Drawbacks of the Envelopment Program.} \subsubsection{Efficiency Frontier of the PPS is 
 not Identified.} Decision variables in the VRS PPS Eq. (\ref{A.2}),   $\lambda_{oj}, \forall j \in$ J, are estimated by the above envelopment program. Hyperplanes Eqs. (\ref{A.4})–(\ref{A.7}) construct the convex feasible space of the envelopment program. The optimal extreme point is located at one of the extreme points of the convex feasible space: the intersection of the objective function and several hyperplanes. On the other hand, DMU-o has a project point on the PPS frontier that is constituted by efficient DMUs. 
 
 \indent The optimal extreme point and projection point belong to the convex feasible space and PPS space. These two points are identified by the decision variables and observed data. The envelopment program cannot identify the projection point because artificial goal weights are still under research. Hence, the estimated intensities and slacks could not be applied to the projection point on the PPS. Similarly, the data envelopment program cannot solve the CRS PPS Eq. (\ref{A.1}).  
\subsubsection{Incomplete Solutions of the Envelopment Program.}
 Each \textit{benchmark condition} Eqs. 
 (\ref{A.4}) and (\ref{A.5}) corresponds to a variable weight, $v_{io}$ (or $u_{ro}$), in the dual (multiplier) program. Using the artificial goal weights Eq. (\ref{A.3}) maximizes the inefficiency score of DMU-o, $F_o^*$. Optimized solutions for envelopment and multiplier programs are used to verify their correctness. As described in the next section, the multiplier program obtains incomplete solutions. Hence, envelopment solutions may be incomplete.  
\subsubsection{The Envelopment Program Does not have Convexity Condition.} Both CRS and VRS envelopment programs have convex feasible spaces. The "convexity condition" Eq. (\ref{A.6}) reduces the space of the CRS envelopment program. It is not equivalent to restricting the non-convex CRS PPS space Eq. (\ref{A.1}) to have the convex VRS PPS space Eq. (\ref{A.2}). Eq. (\ref{A.6}) has a particular scalar 1. In our VGA models, each DMU-j has an intensity $\pi_{oj}$, $\forall j\in J$.  DMU-o managing the SIC scalar $\kappa_o^c$ is significant in the efficiency measurements, where SIC is the condition $\sum_{\forall j \in J} \pi_{oj}=\kappa_o^c$.
 \subsection{Drawbacks of the Multiplier Program.}\label{sec:A2} 
\subsubsection{The Unrealistic Reference Peers.}  The set of artificial goal weights, Eq. (\ref{A.10}), is applied to evaluate every DMU-o,  regardless of its configurations of inputs and outputs. The identified reference peers may not have analogous configurations with DMU-o. The optimized solutions are unrealistic.  
\subsubsection{Inaccurate Solutions of CRS models.}
\indent  
In the CRS model, DMU-o has an efficiency score equal to the ratio between the weighted sum of outputs and the weighted sum of inputs, see Eq. (11) in [Banker et al.]\cite{14},  $\sum_{\forall r\in R} u_{ro} y_{ro}/\sum_{\forall i\in I} v_{io} x_{io}  $. Eq. (\ref{A.8}) aims to minimize DMU-o's excess $f_o^*$, the inefficiency score (=1-efficiency score). The model can produce inaccurate solutions when the estimated virtual input  $\sum_{\forall i\in I} v_{io}^* x_{io}$ is not equal to 1. CRS DEA models use various artificial goal weights and cannot ensure the estimated virtual input equals 1. 
\subsubsection{Inaccurate Solutions of VRS models.} In Eqs. (\ref{A.8}) and (\ref{A.9}) of the VRS model, the inefficiencies scores are affected by $1\sigma_o$, without theoretical support. However, Halická and Trnovská\cite{16} Section 4 concluded that VRS models suffer from duality issues and unbounded solutions. The inefficiency score  $f_o^*$ cannot be obtained.
\subsubsection{Mistakenly Using Dimensionless Measurement} Banker et al.\cite{14} used Figure 2 to depict the particular VRS PPS with one input (x) and one output (y). The line segments passing through points A and E indicate the increasing and decreasing returns to scale by the sign of  $\sigma_o$. Each line respectively intercepts with the vertical and horizontal axis with a length of -$\sigma_o$ and has the same measurement unit of x and y. In this case, Eq. (\ref{A.8}) is expressed as following:
\begin{equation} \label {A.12}
    f_o^\star=\min_{v_o, u_o,\sigma_o} v_o x_o-u_o y_o+1\sigma_o.
\end{equation}
In the above equation, the four elements, $f_o^\star, v_o x_o, u_o y_o$, and $ 1\sigma_o$, should have been measured by the virtual currency \$. Hence, $f_o$ is not the dimensionless inefficiency score. Our VGA models have fixed this problem for cases with multiple inputs and outputs.
\subsubsection{Returns to Scale Is not Measured.}
\indent The actual measurement of the returns to scale in practice relies on the concept of \textit{scale elasticity} and {scale efficiency}, see Page 31 in Sickles and Zelenyuk\cite{5}. However, no one in the DEA literature has presented a solution for multiple inputs and outputs. 
\subsection{Cannot Measure Super-efficiencies.}
\indent The CRS and VRS ARRs DEA models identified Pareto-efficient DMUs with an efficiency score equal to 1. These models are modified to estimate each efficient DMU's super-efficiency, which is more significant than 1. A DMU exhibiting higher super-efficiency surpasses other efficient DMUs. An efficient DMU expands and reduces inputs and outputs, deteriorating its super-efficiency to 1. Our second scenario furnishes super-efficiencies, including and excluding the SIC condition. 

\indent Chen and Du\cite{17} reviewed linear, nonlinear, and integer VRS models, which attempted to measure super-efficiencies. They concluded that convexity constraints caused the infeasibility issue. Unfortunately, those linear programming-based models discussed in [Chen and Du]\cite{17} only present the envelopment program without the multiplier program; those models may not fulfill the duality characteristics. The CRS models have solutions of the envelopment program but without the multiplier program’s solutions. CRS and VRS DEA models have the same problem in measuring inefficiencies and super-efficiencies, even if zero data does not exist.
\subsection{Heterogeneity Caused Unrealistic Evaluations.}
Li et al.\cite{18} raised concerns regarding input configurations significantly differing among DMUs. Aleskerov and Petrushchenko\cite{19} argued that DEA may not be applicable when DMUs operate under vastly different conditions and have significantly distinct input and output compositions. This heterogeneity poses challenges for analysis.

\indent DEA models rely on artificial goal weights and consequently identify reference peers with non-analogous criteria configurations relative to DMU-o. Instead, we propose that DMU-o systematically determines the \textit{goal price} of virtual inputs and outputs to identify reference peers with analogous criteria configurations.
\subsection{Summaries}

This subsection summarizes foundational research and highlights critical limitations for applying DEA models in our study.

\indent Using linear programming frameworks is a notable strength of DEA models, allowing for robust efficiency analyses. However, despite the extensive exploration of artificial goal weights within DEA theories, several critical drawbacks remain:
\begin{romanlist}

\item The solutions provided are incomplete, leaving gaps in the efficiency analysis.
\item DMU-o struggles to identify reference peers with analogous criteria configurations, resulting in unrealistic estimated criteria benchmarks.
\item DMU-o fails to select achievable benchmarks and thus cannot attain a realistic efficiency score.
\item The ranking of Pareto-efficient DMUs according to their super-efficiencies is not feasible.
\end{romanlist}

\indent Modified and extended Additive Models, utilized in various applications, suffer from similar drawbacks and have not been thoroughly discussed in the literature. Researchers such as Dyson et al.\cite{20}, Cooper et. al.\cite{4}, Brown\cite{21}, Zarrin et al.\cite{22}, and Zhu et al.\cite{23} have extensively examined the limitations and methodologies of DEA. Subsequent literature reviews by Mergoni and Witte\cite{24} and Panwar et al.\cite{8} further emphasize the ongoing challenges within DEA models despite numerous refinements. Additionally, Emrouznejad and Yang\cite{25} highlighted the extensive scope of DEA literature, noting the existence of over ten thousand journal articles and a hundred textbooks published in the past fifty years.

These summaries underscore the persistent challenges and the extensive research efforts aimed at addressing the limitations of DEA models, which are crucial for the methodological framework of our study.
\section{Comparisons Between DEA and VGA} 
Table \ref{table:4} summarizes the comparisons of the first scenario, with a footnote providing additional references for further reading. Specifically, DEA models focus on inefficiency scores, while VGA models target virtual gaps. It is noted at the end of the table that DEA's inefficiency scores can be infeasible and unbounded, posing a significant limitation. In contrast, VGA models generate feasible inefficiency scores, effectively addressing this issue. The DMU-o selects the SIC scalar to ensure that these evaluations and benchmarks are practical and attainable.

\indent Charnes et al.\cite{13} and Banker et al.\cite{14} introduced the CRS and VRS DEA models and proposed a two-phase DEA procedure to assess DMU-o. In Phase One, the method obtains overestimated one-sided radial efficiency scores. In Phase Two, the Additive Models measure the remaining inefficiencies. However, the optimal solutions obtained from the two phases cannot be aggregated for further applications.

In contrast, VGA models PT and TSc, which utilize the specific SIC scalar $\kappa_o^c=1$, provide the lower bound efficiency scores of the CRS and VRS DEA models. As illustrated in Appendix A, while CRS and VRS Additive Models offer a robust theoretical basis, they often result in incomplete solutions, highlighting their limitations in practical applications. These models form the foundational frameworks of all DEA models and are essential for understanding efficiency measurements. The first and second VGA models enhance these frameworks by addressing their shortcomings and providing a more detailed and holistic evaluation of efficiency.
\section{ Introduction to Computational Tools}\label{appendix: BBBBB}
In our study, the linear programming (LP) models showcased in Tables \ref{table:2} and \ref{table:3} were solved using Solver, an Add-In for Microsoft Office Excel. This tool facilitates the optimization of models by adjusting the values of the decision variables to meet the constraints and objectives specified by the user. We have meticulously structured the LP models within the accompanying Excel file to enhance understanding and facilitate replication. Each model is comprehensively laid out, showcasing the parameters $(X, Y)$ and $(x_o, y_o)$, decision variables, and constraints. 

\indent The functionality embedded within each solution symbol is detailed, allowing readers to trace the computation logic seamlessly. For clarity, consider the following example from Table 3, DMU-K Evaluation. This evaluation involves models PT-K, TS1-K, and TS2-K. The scalar $\kappa_o^1$, computed in the PT-K model, is a crucial link, hyperlinked to the TS1-K model for further computation. Similarly, the "Sensitivity Report TS1-K" spreadsheet elucidates the calculation of $\kappa_o^2$, which, in turn, is hyperlinked to the TS2-K model.
This interconnected setup not only facilitates a deeper understanding of the computational process but also aids in the transparent communication of how each model's outputs contribute to subsequent calculations.

\indent To ensure the utility of the Excel models, each spreadsheet is designed with readability and ease of navigation in mind. Hyperlinks between related models and their components encourage an intuitive exploration of the computational logic underpinning our research findings. The computations for the second scenario are similar.
 \begin{table}
\centering
\setlength{\tabcolsep}{2pt}
\renewcommand{\arraystretch}{0.91}
  \caption{Comparing DEA vs. VGA Models.}
    \label{table:4} 
\footnotesize
\begin{tabularx}{\linewidth}{|m{0.5cm}|m{1.3cm}|m{2.2cm}|m{2.8cm}|m{2.3cm}|m{2.7cm}|}
\hline
 Row&Compare& \multicolumn{2}{| c |}{Model}&\multicolumn{2}{| c |}{Model}\\
\hline
1&& CRS DEA& First VGA& VRS DEA&Second VGA\\
\hline
2&	Decision variables&	Criteria Weights(1/\dag)&	Criteria virtual unit prices (\$)&	Criteria Weights (1/\dag)&	Criteria virtual unit prices (\$)\\ \hline
3&		&Weighted criteria (-)&	Criteria virtual prices (\$)&	Weighted criteria (-)&Criteria virtual prices  (\$)\\ \hline
4&	Lower bounds&	An artificial goal weight for each weight (1/\dag)&	A unified Systematic goal price for all virtual prices (\$)&	An artificial goal weight for each weight (1/\dag)&	A unified Systematic goal price for all virtual prices (\$)\\ \hline
5&	Decision variables&	Criteria Slacks(\dag)&	Criteria adjustment ratios (-)&	Criteria Slacks(\dag)&	Criteria adjustment ratios (-)\\\hline
6&	Decision variables&	Intensities of DMUs (-)\ddag 1&	Intensities of DMUs (-) \ddag 2 &	Intensities of DMUs (-) \ddag 3&	Intensities of DMUs (-) \ddag 4\\\hline
7&	Decision variable&	&&		SIC Scalar weight (-) \ddag 5 &	SIC Scalar unit price (\$) \ddag 6 \\\hline
&Virtual scalar&&&1$\times$ SIC Scalar weight (-)& Scalar $\times$ SIC Scalar unit price (\$)\\ \hline  
8&	Objective function &	Inefficiency score $F_o^*$ and $f_o^*$ (-)&	Virtual gap (\$)&		Inefficiency score $F_o^*$ and $f_o^*$ (-)&	Virtual gap (\$)\\\hline
9&	Objective values&	 $_o^*$= total of products of slacks and artificial goal weights (-)&	The total adjustment ratios multiply the systematic goal price (\$)&	 $F_o^*$= total of products of slacks and artificial goal weights (-)&	The total adjustment ratios multiply the systematic goal price (\$)\\\hline
10&	Objective values&	 $F_o^*$= virtual input minus output (-)&	Virtual input minus output (\$)&	 $F_o^*$= virtual input minus output plus scalar weight (-)&	Virtual input minus output plus (virtual scalar)= (Affected virtual input)-(Affected virtual Output) (\$) \ddag7\\\hline
11&	0 $<F_o^*\leq 1$&	Not satisfied&	 $F_o^*$ = (Virtual gap)/(Virtual input)	&Unbounded solutions \ddag8&	 $F_o^*$=1-(Affected virtual output) / (Affected virtual input)  \\\hline
12&	Reference peers&	Non-analogous with DMU-o&	Analogous with DMU-o&	Not solved&	Analogous with DMU-o\\\hline
13&	SIC scalar &	Not Available &	Not Available &	=1 \ddag 9 &	Within a systematic range \ddag10\\\hline
14&Scalar selection&No & No& Not allowed& Allowable\ddag11\\\hline
15&	Benchmarks of inputs \& outputs&	Unrealistic& 	Realistic&	Not solved& Realistic \& achievable\\\hline
 \end{tabularx}
Note: Measurement units are indicated in the parenthesize such as (-) dimensionless, (\$) the virtual currency \$, and (\dag) the same as the criterion.
\ddag 1, \ddag 3, and \ddag 5:  see [Eqs. (3.1), (4.38), and (4.43)] in Cooper et al.\cite{4}; \ddag 2, \ddag 4, and \ddag 6 are illustrated in Section \ref{sec:2.1};  Section \ref{sec:2.3} derives \ddag 7.
\ddag8: see Halická and Trnovská\cite{16}. \ddag9 is digested from [Eqs. (4.1) \& (4.37)] of Cooper et al.\cite{4}.
\ddag10 and \ddag11 can be find in Section \ref{sec:4}.\\ \end{table}
\end{document}